# Global Finite-Time Attitude Tracking via Quaternion Feedback

Haichao Gui and George Vukovich*

## Abstract

**This paper addresses the attitude tracking of a rigid body using a quaternion description. Global finite-time attitude controllers are designed with three types of measurements, namely, full states, attitude plus constant-biased angular velocity, and attitude only. In all three scenarios hybrid control techniques are utilized to overcome the well-known topological constraint on the attitude manifold, while coupled nonsmooth feedback inputs are designed via homogeneous theory to achieve finite-time stability. Specially, a finite-time bias observer is derived in the second scenario and a quaternion filter is constructed to provide damping in the absence of velocity feedback. The proposed methods ensure bounded control torques *a priori* and, in particular, include several existing attitude controllers as special cases.**

## 1 Introduction

Attitude control of a rigid body is a fundamental control problem that has been extensively studied. The attitude configuration of a rigid body SO(3), consisting of 3×3 rotation matrices, is a compact noncontractible manifold. This distinctive feature precludes the existence of globally stabilizing continuous state-feedback laws [1], making global attitude control a challenging issue.

Unit quaternions, as pervasive, global, nonsingular attitude coordinates, cover SO(3) twice. As a result, some quaternion-based attitude control systems can give rise to two antipodal equilibria representing the same desired attitude. If one equilibrium is stable while the other is not, the unwinding phenomenon can occur, yielding an unnecessary full rotation even for small initial attitude errors (e.g., the methods in [1-3]). The approaches to resolve this problem can be categorized into two types. The first type allows continuous state-feedback laws a quaternion manifold or SO(3) but entails unstable equilibria other than the desired attitude, thus achieving almost global stabilization [4]. The other relies on discontinuous switching laws to determine a direction to reach the desired attitude. Examples of this type can be found in [5,6], where a memoryless switching law was utilized, and in [7-10], where a hysteresis mechanism was constructed in the framework of hybrid system theory to further mitigate the chattering due to perturbations such as measurement noise.

All the above attitude control schemes at best produce exponential convergence with infinite settling time, as opposed to finite-time control with finite convergence time. Finite-time stable systems usually demonstrate fast convergence rates and significant disturbance rejection properties [11-14]. Given these properties, finite-time attitude regulators were constructed in [15] via a fractional-power feedback domination approach, in [16] by the terminal sliding mode (TSM) method, in [17,18] by homogeneous theory, and in [19] by a TSM-like method proposed in [20] for simple mechanical systems. None of these control laws, however, ensures global stability. Although a shorter rotation path is prespecified in [15] according to the initial attitude, this method still fails to produce a well-defined vector field on $SO(3)\times\mathbb{R}^3$ [1] after initiation and can exhibit unwinding due to neglecting the effect of the initial kinetic energy (or, velocity), as shown in [7]. The controllers in [18,19], respectively designed on quaternion manifold and SO(3), both ensure almost global finite-time stability (AGFTS) with continuous inputs. They avoid unwinding but introduce unstable equilibria on the manifold of $\pi$ rotations. As a result, the response can be sluggish if the initial condition is close to or coincides with one of these nontrivial equilibria.

In this paper, global finite-time attitude controllers (GFTACs) for tracking maneuvers are proposed based on quaternion representations for three measurement scenarios: full states, attitude plus angular velocity corrupted by an unknown constant bias, and attitude alone. In all three cases finite-time feedback inputs are constructed by a homogenous method and injected in a manner such that the closed-loop system involves only two antipodal equilibria on a quaternion manifold. As a result, simple hysteresis-based switching laws can be further integrated through hybrid control techniques to enable global stability of the two disconnected quaternion equilibria. In addition, the design incorporates a finite-time observer to estimate the velocity bias for the second scenario and a quaternion filter to provide damping for the velocity-free case. The resulting control laws possess a simple, nonlinear, proportional-derivative (PD) structure, ensure bounded control torques *a priori*, and in particular, can recover the methods in [2,3,5-7] by proper selection of control parameters. Our methods avoid the unwinding problem, as also guaranteed by the control algorithms developed on SO(3) [4,9,10,19,21], and achieves global finite-time convergence, which is not a feature of the prior controls in [4,9,10,19,21].

* H. Gui and G. Vukovich are with the Department of Earth and Space Science and Engineering, York University, Toronto, Ontario M3J 1P3, Canada. (e-mail: guihaichao@gmail.com; vukovich@yorku.ca; tel: +1-416-736-2100x30090).

## 2 Preliminaries and System Models

### 2.1 Definitions and Lemmas

For any $x \in \mathbb{R}$ and $\alpha \ge 0$, let $\mathrm{sgn}^{\alpha}(x) = \mathrm{sgn}(x)|x|^{\alpha}$ and $\mathrm{sat}_{\alpha}(x) = \mathrm{sgn}(x)\min\{|x|^{\alpha}, 1\}$, where $\mathrm{sgn}(\cdot)$ is the standard sign function. Clearly, $\mathrm{sgn}^{\alpha}(x)$ is a continuous nonsmooth function if $0 < \alpha < 1$, while $\mathrm{sat}_{\alpha}(x)$ becomes the standard saturation function $\mathrm{sat}(x)$ if $\alpha = 1$. For $\forall \boldsymbol{x} \in \mathbb{R}^n$, let $\mathrm{sgn}^{\alpha}(\boldsymbol{x}) = [\mathrm{sgn}^{\alpha}(x_1), \cdots, \mathrm{sgn}^{\alpha}(x_n)]^T$ and $\mathrm{sat}_{\alpha}(\boldsymbol{x}) = [\mathrm{sat}_{\alpha}(x_1), \cdots, \mathrm{sat}_{\alpha}(x_n)]^T$. In addition, $y = \mathcal{O}(x)$ means $|y| \le c|x|$ for sufficiently small $|x|$ and some constant $c > 0$ while $y = o(x)$ means $\lim_{x\to 0} y/x = 0$. Denote by $\|\cdot\|$ the Euclidean norm and $\mathbb{I}_n$ the index set $\{1, \cdots, n\}$. Given $\varepsilon > 0$ and a weight vector $\boldsymbol{r} = (r_1, \cdots, r_n)$, where $r_i > 0$, $\forall i \in \mathbb{I}_n$, define a dilation operator $\Delta_{\varepsilon}^{r}$ as $\Delta_{\varepsilon}^{r}\boldsymbol{x} = [\varepsilon^{r_1}x_1, \cdots, \varepsilon^{r_n}x_n]^T$ for $\forall \boldsymbol{x} \in \mathbb{R}^n$ [22]. To deal with time-dependent functions and systems, the dilation operator $\Delta_{\varepsilon}^{r}$ is extended as $\Delta_{\varepsilon}^{r}(\boldsymbol{x}, t) = (\Delta_{\varepsilon}^{r}\boldsymbol{x}, t)$. Given a function $V(\boldsymbol{x}): \mathbb{R}^n \mapsto \mathbb{R}$ and a vector field $\boldsymbol{f}(\boldsymbol{x}): \mathbb{R}^n \mapsto \mathbb{R}^n$, denote $L_{\boldsymbol{f}}V(\boldsymbol{x})$ as the Lie derivative of $V$ along $\boldsymbol{f}(\boldsymbol{x})$. The function $V(\boldsymbol{x})$ is said to be homogeneous of degree $k \in \mathbb{R}$ with respect to $\Delta_{\varepsilon}^{r}$ if $V(\Delta_{\varepsilon}^{r}\boldsymbol{x}) = \varepsilon^{k}V(\boldsymbol{x})$.

**Definition 2.1** [23]. Consider the system

$$\dot{\boldsymbol{x}} = \boldsymbol{f}(\boldsymbol{x}, t),\ \boldsymbol{f}(0, t) = 0,\ \boldsymbol{x} \in \mathbb{R}^n, \tag{1}$$

where $\boldsymbol{f}(\boldsymbol{x}, t) = [f_1(\boldsymbol{x}, t), \cdots, f_n(\boldsymbol{x}, t)]^T : U \times \mathbb{R}_{\ge 0} \mapsto \mathbb{R}^n$ is continuous on an open neighborhood $U$ of the origin. Then, $\boldsymbol{f}(\boldsymbol{x}, t)$ is said to be homogeneous of degree $k \in \mathbb{R}$ with respect to a dilation $\Delta_{\varepsilon}^{r}$ if $f_i(\Delta_{\varepsilon}^{r}\boldsymbol{x}, t) = \varepsilon^{r_i + k} f_i(\boldsymbol{x}, t)$ for $\forall i \in \mathbb{I}_n$, $\forall \boldsymbol{x} \in U$, and any $\varepsilon > 0$. System (1) is said to be homogeneous if $\boldsymbol{f}(\boldsymbol{x}, t)$ is homogeneous.

**Definition 2.2** [14]. Consider system (1) and denote by $U$ a neighborhood of $\boldsymbol{x} = 0$. The origin is uniformly finite-time stable if it is 1) uniformly Lyapunov stable in $U$ and 2) uniformly finite-time convergent in $U$ (that is, there exists $T(\boldsymbol{x}_0) \ge 0$, depending only on $\boldsymbol{x}_0 \in U$, such that the system solution denoted by $\boldsymbol{x}(t, \boldsymbol{x}_0, t_0)$ satisfies $\lim_{t \to T(\boldsymbol{x}_0)} \boldsymbol{x}(t, \boldsymbol{x}_0, t_0) = 0$ and $\boldsymbol{x}(t, \boldsymbol{x}_0, t_0) = 0$ for $t \ge T(\boldsymbol{x}_0) + t_0$). If $U = \mathbb{R}^n$, then the origin is uniformly globally finite-time stable (UGFTS). If system (1) is time-invariant, the qualifier '*uniformly*' can be omitted in the preceding statements.

The following result is an extension of Lemma 3 in [12], dealing with autonomous systems, to non-autonomous systems. The condition that (3) uniformly holds for $t \ge 0$ enables us to prove Lemma 2.1 in a manner analogous to that of Theorem 3 in [22].

**Lemma 2.1**. Consider the system

$$\dot{\boldsymbol{x}} = \boldsymbol{f}(\boldsymbol{x}) + \hat{\boldsymbol{f}}(\boldsymbol{x}, t),\ \boldsymbol{f}(0) = 0,\ \boldsymbol{x} \in \mathbb{R}^n, \tag{2}$$

where $\boldsymbol{f}(\boldsymbol{x})$ is a continuous homogeneous vector field of degree $k < 0$ with respect to a dilation $\Delta_{\varepsilon}^{r}$, and $\hat{\boldsymbol{f}}(\boldsymbol{x}, t)$ satisfies $\hat{\boldsymbol{f}}(0, t) = 0$. Assume that $\boldsymbol{x} = 0$ is an asymptotically stable equilibrium of the system $\dot{\boldsymbol{x}} = \boldsymbol{f}(\boldsymbol{x})$. Then, $\boldsymbol{x} = 0$ is a uniformly locally finite-time stable equilibrium of system (2) if

$$\lim_{\varepsilon \to 0} \frac{\hat{f}_i(\Delta_{\varepsilon}^{r}\boldsymbol{x}, t)}{\varepsilon^{r_i + k}} = 0,\ i \in \mathbb{I}_n,\ \forall \boldsymbol{x} \ne 0, \tag{3}$$

holds uniformly for $t \ge 0$. Moreover, if system (2) is uniformly globally asymptotically stable (UGAS) and uniformly locally finite-time stable, it is UGFTS.

**Proof**. See Appendix A. ■

**Remark 2.1**. As shown in Appendix A, the conditions of Lemma 2.1 ensure the existence of a $C^1$ locally strict Lyapunov function $V(\boldsymbol{x}) \ge 0$ for system (2) such that $\dot{V}(\boldsymbol{x}, t) \le -c_1 V^{\beta}(\boldsymbol{x})$ holds for all $(\boldsymbol{x}, t) \in U \times \mathbb{R}_{\ge 0}$, where $c_1 > 0$ is a constant, $U$ is a neighborhood of the origin, and $0 < \beta < 1$ is determined by $k$, the homogeneous degree of $\boldsymbol{f}(\boldsymbol{x})$. Now, assume that system (2) is perturbed to $\dot{\boldsymbol{x}} = \boldsymbol{f}(\boldsymbol{x}) + \hat{\boldsymbol{f}}(\boldsymbol{x}, t) + \boldsymbol{d}(t)$ by a small perturbation $\boldsymbol{d}(t)$ satisfying $|L_{\boldsymbol{d}}V(\boldsymbol{x}, t)| \le \delta < c_1$ for $(\boldsymbol{x}, t) \in U \times \mathbb{R}_{\ge 0}$. Noting that $\dot{V} = L_{\boldsymbol{f}}V + L_{\hat{\boldsymbol{f}}}V + L_{\boldsymbol{d}}V$, it follows that $\dot{V}(\boldsymbol{x}, t) \le -c_1 V^{\beta}(\boldsymbol{x}) + \delta$ and thus $\boldsymbol{x}$ is stabilized to the region $V(\boldsymbol{x}) \le (\delta/c_1)^{1/\beta}$. Since $0 < \delta/c_1 < 1$ and $1/\beta > 1$, the ultimate bound on the state can be reduced to a higher order than the bound on the perturbation by properly tuning $k$ to decrease $\beta$. If system (2) is exponentially stable, the size of the convergent region is $\delta/c_1$, only the same order as the perturbation magnitude. Therefore, it can be expected that the proposed finite-time methods, which have an adjustable homogeneous degree, achieve better robustness than the asymptotic or exponential control algorithms in [2-10,21] against small perturbations imposed on the system

### 2.2 Equations of Attitude Motion

Denote by $\boldsymbol{Q} = [q_0, \boldsymbol{q}^T]^T \in \mathbb{R}^4$ a quaternion, where $q_0 \in \mathbb{R}$ and $\boldsymbol{q} = [q_1, q_2, q_3]^T \in \mathbb{R}^3$ are usually called the scalar and vector parts, respectively. The multiplication of two quaternions $\boldsymbol{Q} = [q_0, \boldsymbol{q}^T]^T$ and $\boldsymbol{P} = [p_0, \boldsymbol{p}^T]^T$ is defined as

$$\boldsymbol{Q}\otimes\boldsymbol{P}=\begin{bmatrix} q_0p_0-\boldsymbol{q}^T\boldsymbol{p} \\ q_0\boldsymbol{p}+p_0\boldsymbol{q}+\boldsymbol{q}\times\boldsymbol{p}\end{bmatrix}, \tag{4}$$

where $\times$ is the cross product in $\mathbb{R}^3$. The quaternion multiplication is associative and distributive but is not commutative. The conjugation of a quaternion $\boldsymbol{Q}\in\mathbb{R}^4$ is defined as $\boldsymbol{Q}^*=[q_0,-\boldsymbol{q}^T]^T\in\mathbb{R}^4$. Note that $(\boldsymbol{Q}\otimes\boldsymbol{P})^*=\boldsymbol{P}^*\otimes\boldsymbol{Q}^*$. With the identity element $\mathbf{1}=[1,0,0,0]^T$, a unit quaternion can be defined as $\boldsymbol{Q}\in\mathbb{S}^3=\{\boldsymbol{Q}\in\mathbb{R}^4:\boldsymbol{Q}\otimes\boldsymbol{Q}^*=\mathbf{1}\}$. Alternatively, $\boldsymbol{Q}=[q_0,\boldsymbol{q}^T]^T\in\mathbb{S}^3$ can also be expressed as $q_0=\cos(\phi/2)$ and $\boldsymbol{q}=\sin(\phi/2)\boldsymbol{\eta}$, where $\phi\in[0,2\pi]$ and $\boldsymbol{\eta}$ is a unit vector, also called as the eigenaxis vector.

Letting $\boldsymbol{Q}\in\mathbb{S}^3$ and $0\le\alpha<1$, two functions of $\boldsymbol{Q}$ are defined as

$$\boldsymbol{\kappa}_0(\boldsymbol{Q},\alpha)=\begin{cases}\dfrac{\boldsymbol{q}}{\|\boldsymbol{q}\|^{\alpha}}, & \|\boldsymbol{q}\|\neq 0\\ 0\in\mathbb{R}^3, & \|\boldsymbol{q}\|=0\end{cases},\quad \boldsymbol{\kappa}_1(\boldsymbol{Q},\alpha)=\begin{cases}\dfrac{\boldsymbol{q}}{\left(\sqrt{2(1-q_0)}\right)^{\alpha}}, & q_0\neq 1\\ 0\in\mathbb{R}^3, & q_0=1\end{cases}. \tag{5}$$

The subtraction of $\boldsymbol{\kappa}_1$ and $\boldsymbol{\kappa}_0$ is denoted by $\bar{\boldsymbol{\kappa}}(\boldsymbol{Q},\alpha)=\boldsymbol{\kappa}_1(\boldsymbol{Q},\alpha)-\boldsymbol{\kappa}_0(\boldsymbol{Q},\alpha)$. These functions possess the following important properties whose proofs can be found in Appendix B.

**Property 1**. $\boldsymbol{\kappa}_0(\boldsymbol{Q},\alpha)$ and $\boldsymbol{\kappa}_1(\boldsymbol{Q},\alpha)$ are continuous for every $\boldsymbol{Q}\in\mathbb{S}^3$ and $\|\boldsymbol{\kappa}_1(\boldsymbol{Q},\alpha)\|\le 1$.

**Property 2**. If $\|\boldsymbol{q}\|\ll 1$ and $q_0>0$, it follows that $\bar{\boldsymbol{\kappa}}(\boldsymbol{Q},\alpha)\approx-\alpha\|\boldsymbol{q}\|^2\boldsymbol{\kappa}_0(\boldsymbol{Q},\alpha)/8$.

Let $\mathcal{F}_{\mathcal{I}}$ represent the inertial frame and $\mathcal{F}_{\mathcal{B}}$ denote the body-fixed frame of a rigid body. The attitude of $\mathcal{F}_{\mathcal{B}}$ relative to $\mathcal{F}_{\mathcal{I}}$ can be represented by a unit quaternion $\boldsymbol{Q}\in\mathbb{S}^3$ and the rotation matrix from $\mathcal{F}_{\mathcal{I}}$ to $\mathcal{F}_{\mathcal{B}}$ can be expressed as $\boldsymbol{R}(\boldsymbol{Q})=(q_0^2-\boldsymbol{q}^T\boldsymbol{q})\boldsymbol{I}_3+2\boldsymbol{q}\boldsymbol{q}^T-2q_0\boldsymbol{q}^{\times}$, where $\boldsymbol{I}_3$ denotes the 3×3 identity matrix and the operator $\boldsymbol{x}^{\times}$ denotes a skew-symmetric matrix satisfying $\boldsymbol{x}^{\times}\boldsymbol{y}=\boldsymbol{x}\times\boldsymbol{y}$ for $\boldsymbol{x},\boldsymbol{y}\in\mathbb{R}^3$. Letting $\boldsymbol{x}\in\mathbb{R}^3$ be a vector expressed in $\mathcal{F}_{\mathcal{I}}$, it can be verified that $\boldsymbol{R}(\boldsymbol{Q})\boldsymbol{x}=\boldsymbol{Q}^*\otimes\boldsymbol{x}\otimes\boldsymbol{Q}$.

The attitude kinematics and dynamics of a rigid body in terms of unit quaternions are given by

$$\dot{\boldsymbol{Q}}=\frac{1}{2}\boldsymbol{Q}\otimes\boldsymbol{\omega}=\frac{1}{2}\begin{bmatrix}-\boldsymbol{q}^T\\ \boldsymbol{E}(\boldsymbol{q})\end{bmatrix}\boldsymbol{\omega};\ \boldsymbol{E}(\boldsymbol{q})\triangleq\boldsymbol{q}^{\times}+q_0\boldsymbol{I}_3, \tag{6}$$

$$\boldsymbol{J}\dot{\boldsymbol{\omega}}=-\boldsymbol{\omega}\times\boldsymbol{J}\boldsymbol{\omega}+\boldsymbol{u}, \tag{7}$$

where $\boldsymbol{\omega}\in\mathbb{R}^3$ is the angular velocity with respect to $\mathcal{F}_{\mathcal{I}}$ expressed in $\mathcal{F}_{\mathcal{B}}$, $\boldsymbol{J}\in\mathbb{R}^{3\times3}$ is the symmetric rigid-body inertia matrix, and $\boldsymbol{u}\in\mathbb{R}^3$ represents the control torque.

### 2.3 Equations of Relative Attitude Motion

The desired attitude motion is generated by a frame $\mathcal{F}_{\mathcal{D}}$ with its attitude and angular velocity denoted by $\boldsymbol{Q}_d=[q_{d0},\boldsymbol{q}_d^T]^T\in\mathbb{S}^4$ and $\boldsymbol{\omega}_d\in\mathbb{R}^3$, which satisfy the kinematics given in (6). The error quaternion and error angular velocity of $\mathcal{F}_{\mathcal{B}}$ with respect to $\mathcal{F}_{\mathcal{D}}$ can then be defined as

$$\boldsymbol{Q}_e=\boldsymbol{Q}_d^*\otimes\boldsymbol{Q}, \tag{8}$$

$$\boldsymbol{\omega}_e=\boldsymbol{\omega}-\boldsymbol{R}(\boldsymbol{Q}_e)\boldsymbol{\omega}_d=\boldsymbol{\omega}-\bar{\boldsymbol{\omega}}_d;\ \bar{\boldsymbol{\omega}}_d\triangleq\boldsymbol{R}(\boldsymbol{Q}_e)\boldsymbol{\omega}_d, \tag{9}$$

where $\boldsymbol{R}(\boldsymbol{Q}_e)$ is the transformation matrix from $\mathcal{F}_{\mathcal{D}}$ to $\mathcal{F}_{\mathcal{B}}$. With these definitions of $\boldsymbol{Q}_e$ and $\boldsymbol{\omega}_e$, the equations of attitude motion relative to $\mathcal{F}_{\mathcal{D}}$ can be written as

$$\dot{\boldsymbol{Q}}_e=\frac{1}{2}\boldsymbol{Q}_e\otimes\boldsymbol{\omega}_e, \tag{10}$$

$$\boldsymbol{J}\dot{\boldsymbol{\omega}}_e=\Xi(\boldsymbol{\omega}_e,\bar{\boldsymbol{\omega}}_d)\boldsymbol{\omega}_e-\bar{\boldsymbol{\omega}}_d^{\times}\boldsymbol{J}\bar{\boldsymbol{\omega}}_d-\boldsymbol{J}\boldsymbol{R}(\boldsymbol{Q}_e)\dot{\boldsymbol{\omega}}_d+\boldsymbol{u}, \tag{11}$$

where $\Xi(\boldsymbol{\omega}_e,\bar{\boldsymbol{\omega}}_d)=(\boldsymbol{J}(\boldsymbol{\omega}_e+\bar{\boldsymbol{\omega}}_d))^{\times}-\bar{\boldsymbol{\omega}}_d^{\times}\boldsymbol{J}-\boldsymbol{J}\bar{\boldsymbol{\omega}}_d^{\times}$ is skew-symmetric. In the following, $\boldsymbol{u}_d=\bar{\boldsymbol{\omega}}_d^{\times}\boldsymbol{J}\bar{\boldsymbol{\omega}}_d+\boldsymbol{J}\boldsymbol{R}(\boldsymbol{Q}_e)\dot{\boldsymbol{\omega}}_d$, which represents the torque required to maintain the desired attitude motion with zero tracking error.

**Assumption 2.1**. The desired angular velocity and angular acceleration are continuous and bounded as $\|\boldsymbol{\omega}_d(t)\|\le\varpi_1$ and $\|\dot{\boldsymbol{\omega}}_d(t)\|\le\varpi_2$, where $\varpi_1$ and $\varpi_2$ are known constants.

## 3 Main Results

In this section, control torques are designed such that $(\boldsymbol{Q}_e, \boldsymbol{\omega}_e) \to (\pm\mathbf{1}, 0)$ in finite time. Three cases with full-state measurements, measurements of attitude plus angular velocity with unknown constant bias, and attitude-only measurements respectively are considered.

Similarly to [7], hysteretic switching laws, embodied through a binary logic variable $h \in \mathbb{H} \triangleq \{-1,1\}$ as well as a constant hysteresis gap $\delta \in (0,1)$, are used to determine the desired quaternion equilibrium from $\pm\mathbf{1}$. In addition, define an outer semicontinuous set-valued map $\overline{\mathrm{sgn}}(x) \in \mathbb{H}$ for $x \in \mathbb{R}$, where $\overline{\mathrm{sgn}}(x) = \mathrm{sgn}(x)$ for $x \neq 0$ and $\overline{\mathrm{sgn}}(0) \in \mathbb{H}$ [7]. To account for discontinuous switches, the design in the following is posed in the setting of hybrid systems, which are combinations of continuous dynamics, defined on a flow set, and discrete dynamics, defined on a jump set. Following the framework of [24], a hybrid system is denoted by

$$\mathcal{H} \triangleq \begin{cases} \dot{\boldsymbol{x}} = F(\boldsymbol{x}), & \boldsymbol{x} \in C \\ \boldsymbol{x}^+ = G(\boldsymbol{x}), & \boldsymbol{x} \in D \end{cases}, \tag{12}$$

where the flow map $F: \mathbb{R}^n \rightrightarrows \mathbb{R}^n$ governs continuous evolution of the state $\boldsymbol{x}$ when $\boldsymbol{x}$ falls on the flow set $C$, while the jump map $G: \mathbb{R}^n \rightrightarrows \mathbb{R}^n$ governs the discrete dynamics over the jump set $D$ and $\boldsymbol{x}^+$ denotes the state value immediately after a jump.

### 3.1 GFTAC with Full-State Measurements

Here it is assumed that both $\boldsymbol{Q}$ and $\boldsymbol{\omega}$ are available. Let $\boldsymbol{x}_1 = (\boldsymbol{Q}_e, \boldsymbol{\omega}_e, h) \in M \triangleq \mathbb{S}^3 \times \mathbb{R}^n \times \mathbb{H}$. The objective is to globally stabilize the set $E_1 = \{\boldsymbol{x}_1 \in M : \boldsymbol{Q}_e = h\mathbf{1}, \boldsymbol{\omega}_e = 0\}$ in finite time. Inspired by [7], the flow and jump sets are defined respectively as

$$C_1 = \{\boldsymbol{x}_1 \in M : hq_{e0} \geq -\delta\}, \tag{13}$$

$$D_1 = \{\boldsymbol{x}_1 \in M : hq_{e0} \leq -\delta\}, \tag{14}$$

where $C_1 \cup D_1 = M$. A hybrid attitude tracking controller is then designed as

$$\boldsymbol{u}(\boldsymbol{x}_1) = \boldsymbol{u}_d - k_1 \boldsymbol{\kappa}_1(h\boldsymbol{Q}_e, 1-\alpha_1) - k_2\, \mathrm{sat}_{\alpha_2}(\boldsymbol{\omega}_e),\ \boldsymbol{x}_1 \in C_1, \tag{15}$$

$$\boldsymbol{x}_1^+ = G_1(\boldsymbol{x}_1) = (\boldsymbol{Q}_e, \boldsymbol{\omega}_e, \overline{\mathrm{sgn}}(q_{e0})),\ \boldsymbol{x}_1 \in D_1, \tag{16}$$

where $k_1, k_2 > 0$, $0 < \alpha_1 < 1$, and $\alpha_2 = 2\alpha_1/(1+\alpha_1)$. Note that $h$ reverses its sign only on the jump set $D_1$; otherwise, it remains constant and thus $\dot{h} = 0$.

Actually, (16) presents a switching law to change the desired quaternion equilibrium (and thus the desired rotation direction) when the amount of sign mismatch between $h$ and $q_{e0}$ reaches a prespecified hysteresis width $\delta$. In contrast to possible discontinuous switches of $h$, $(\boldsymbol{Q}_e, \boldsymbol{\omega}_e)$ evolves continuously on the entire state space. As can be observed from (15), the control torque consists of a nonlinear, nonsmooth, PD feedback of error states, which ensures closed-loop stability, and a feedforward compensation for the acceleration of the desired trajectory. In addition, the saturation function in (15) imposes an upper bound on the angular velocity feedback term.

**Theorem 3.1**. Consider the hybrid control system given by (10), (11), and (13)-(16) with $k_1, k_2 > 0$, $0 < \alpha_1 < 1$, $\alpha_2 = 2\alpha_1/(1+\alpha_1)$, and $\delta \in (0,1)$. Then, the compact set $E_1$ is UGFTS.

**Remark 3.1**. Following Assumption 2.1, $\|\bar{\boldsymbol{\omega}}_d^\times \boldsymbol{J}\bar{\boldsymbol{\omega}}_d + \boldsymbol{J}\boldsymbol{R}(\boldsymbol{Q}_e)\dot{\boldsymbol{\omega}}_d\| \leq (\varpi_1^2 + \varpi_2)\|\boldsymbol{J}\|$ since the desired trajectory is assumed to be bounded. Property 1 can be employed to show that the components of $\boldsymbol{u}(\boldsymbol{x}_1)$ are bounded by $|u_i(\boldsymbol{x}_1)| < k_1 + k_2 + (\varpi_1^2 + \varpi_2)\|\boldsymbol{J}\|$, $i \in \mathbb{I}_3$. This feature can facilitate the accommodation of actuator saturation constraints *a priori*. On the other hand, the preceding full-state tracking law becomes the hybrid control law in [7] by setting $\alpha_1 = \alpha_2 = 1$ and replacing $\mathrm{sat}_{\alpha_2}(\boldsymbol{\omega}_e)$ with a strongly passive function $\Phi(\boldsymbol{\omega}_e)$, and the smooth controller in [2] by further setting $h \equiv 1$ and omitting the saturation function. The methods in [2,7], however, can only achieve asymptotic or local exponential convergence of the attitude tracking error. In the attitude regulation case, controller (15) is completely independent of the inertia matrix since $\boldsymbol{u}_d \equiv 0$.

**Remark 3.2**. Theoretically, it is possible to extend the almost global finite-time attitude controllers in [18,19] into GFTACs by constructing a synergistic potential function (SPF) via the angular wrapping method [25] for each undesired equilibrium and then hysteretically switching between these SPFs [9]. The GFTAC given by (15) and (16), however, is not such a direct extension of [18]. Note that the quaternion-based feedback in [18] introduces 26 undesired equilibria, which necessitate 26 SPFs and thus imply great complexity in applying the synergistic hybrid technique. The method in [19] involves merely three undesired equilibria on SO(3) but its extension to a GFTAC is still difficult due to the complexity of the construction of SPFs on SO(3) and the determination of the synergistic gap, as shown in [9]. Recently, Lee [10] proposed novel error functions on SO(3) to obtain an explicit synergistic gap and a switching logic with three control modes. In contrast, by means of a carefully designed nonsmooth function $\boldsymbol{\kappa}_1(\cdot,\cdot)$ in (5) controller (15) injects feedback in a manner such that it produces only two antipodal equilibria ($E_1$) representing the same desired attitude. It can be seen by setting $h \equiv 1$ or $h \equiv -1$ (leading to continuous inputs) that one of them is uniformly locally finite-time stable while the other is unstable. As a result, the simple switching logic in (16) can be incorporated to globally stabilize $E_1$ in finite time and only two control modes are involved.

**Proof of Theorem 3.1**. For $\alpha \ge 0$ and $|x| \le 1$, define two functions as $\varphi(x,\alpha) = \left(\sqrt{2(1-x)}\right)^{\alpha}$ and $\rho(x,\alpha) = \varphi(|x|,\alpha) - \varphi(x,\alpha)$. A Lyapunov candidate function is then constructed as

$$V_1(\boldsymbol{x}_1) = \frac{1}{2}\boldsymbol{\omega}_e^T \boldsymbol{J}\boldsymbol{\omega}_e + \frac{2k_1}{1+\alpha_1}\varphi(hq_{e0}, 1+\alpha_1), \tag{17}$$

which satisfies $V_1(\boldsymbol{x}_1) \ge 0$ and $V_1(\boldsymbol{x}_1) = 0$ if and only if $\boldsymbol{x}_1 \in E_1$. To show that the closed-loop system satisfies the sufficient conditions of Lemma 2.1, the following proof is divided into three steps.

*Step 1) Uniform Lyapunov Stability*. Straight computation produces

$$\frac{\mathrm{d}}{\mathrm{d}t}\varphi(hq_{e0}, 1+\alpha_1) = \frac{1+\alpha_1}{2}\boldsymbol{\omega}_e^T \boldsymbol{\kappa}_1(h\boldsymbol{Q}_e, 1-\alpha_1). \tag{18}$$

By means of (18) and noting $\boldsymbol{\omega}_e^T \Xi(\boldsymbol{\omega}_e, \bar{\boldsymbol{\omega}}_d)\boldsymbol{\omega}_e = 0$, the time derivative of $V_1$ along the flow dynamics can then be computed as

$$\dot{V}_1 = -k_2 \boldsymbol{\omega}_e^T \mathrm{sat}_{\alpha_2}(\boldsymbol{\omega}_e) \le 0, \ \boldsymbol{x}_1 \in C_1. \tag{19}$$

Over a jump one can obtain

$$V_1(G(\boldsymbol{x}_1)) - V_1(\boldsymbol{x}_1) = \frac{2k_1}{1+\alpha_1}\rho(hq_{e0}, 1+\alpha_1).$$

Note that $\rho(x,\alpha)$ is monotonically increasing with $x$ when $-1 \le x \le 0$. Since $hq_{e0} \le -\delta$, it follows that $\rho(hq_{e0},\cdot) \le \rho(-\delta,\cdot)$ and hence $V_1(G_1(\boldsymbol{x}_1)) - V_1(\boldsymbol{x}_1) \le -\sigma_1 < 0$, where $\sigma_1 = -2k_1\rho(-\delta, 1+\alpha_1)/(1+\alpha_1) > 0$. The above analysis shows that $V_1(\boldsymbol{x}_1)$ monotonically decreases along the system flows and strictly decreases over jumps. Therefore, $E_1$ is uniformly Lyapunov stable.

*Step 2) Uniform Asymptotic Convergence*. Let $t_j > 0$ denote a time, before which $j \in \mathbb{N}$ jumps have occurred. The analysis in Step 1) indicates $V_1(\boldsymbol{x}_1(t_j)) \le V_1(\boldsymbol{x}_1(0)) - j\sigma_1$. Since $V_1(\boldsymbol{x}_1) \ge 0$, it follows that $j \le V_1(\boldsymbol{x}_1(0))/\sigma_1$, implying that the number of jumps is finite for any bounded initial condition. Thus, the closed-loop trajectory remains in the flow set $C_1$ for some finite time and we can simply assume $\boldsymbol{x}_1 \in C_1$ in the following proof. Equation (19) implies that $\boldsymbol{\omega}_e(t)$ is uniformly bounded and $\lim_{t\to\infty} V_1(t) = V_1(\infty)$ exists and is finite. In addition, $\dot{\boldsymbol{\omega}}_e(t)$ is uniformly bounded and, as a result, $\boldsymbol{\omega}_e(t)$ and thus $\dot{V}_1(t)$ are uniformly continuous. Barbalat's lemma [26] can then be used to conclude that $\dot{V}_1(t) \to 0$ and thus $\boldsymbol{\omega}_e(t) \to 0$. Next, note that $\boldsymbol{\kappa}_1(h\boldsymbol{Q}_e, 1-\alpha_1)$ is a continuous function of $\boldsymbol{Q}_e$ by Property 1 and $\mathrm{sat}_{\alpha_2}(\boldsymbol{\omega}_e)$ is a continuous function of $\boldsymbol{\omega}_e(t)$. On the other hand, $\boldsymbol{Q}_e(t)$ and $\boldsymbol{\omega}_e(t)$ are both uniformly continuous functions of $t$ and uniformly bounded. Therefore, $\boldsymbol{\kappa}_1(h\boldsymbol{Q}_e(t), 1-\alpha_1)$ and $\mathrm{sat}_{\alpha_2}(\boldsymbol{\omega}_e(t))$ are uniformly continuous with respect to $t$ as well. This conclusion stems from the fact that a continuous function over a compact space is uniformly continuous. Substituting (15) into (11), it follows that $\dot{\boldsymbol{\omega}}_e(t)$ is uniformly continuous since all terms involved are uniformly continuous. Again, by Barbalat's lemma, one can deduce that $\dot{\boldsymbol{\omega}}_e(t) \to 0$ and thus $\boldsymbol{\kappa}_1(h\boldsymbol{Q}_e(t), 1-\alpha_1) \to 0$, which implies $\boldsymbol{Q}_e(t) \to h\mathbf{1}$ since $hq_{e0} \ge -\delta$. Recalling the result in Step 1), one obtains that $E_1$ is UGAS.

*Step 3) Uniform Local Finite-Time Stability*. Since the number of jumps is uniformly bounded and the closed-loop trajectory asymptotically converges to $E_1$. We can restrict the analysis below to within a neighborhood of $E_1$ such that $|\omega_{ei}| \le 1$, $i \in \mathbb{I}_3$ and $hq_{e0} > 0$. The closed-loop flow dynamics can then be written as

$$\dot{\boldsymbol{q}}_e = 0.5h\boldsymbol{\omega}_e + \hat{\boldsymbol{f}}_1(\boldsymbol{x}), \tag{20}$$

$$\dot{\boldsymbol{\omega}}_e = -\boldsymbol{J}^{-1}[k_1\boldsymbol{\kappa}_0(h\boldsymbol{Q}_e, 1-\alpha_1) + k_2\,\mathrm{sgn}^{\alpha_2}(\boldsymbol{\omega}_e)] + \hat{\boldsymbol{f}}_2(\boldsymbol{x}, t), \tag{21}$$

where $\boldsymbol{x} = (\boldsymbol{q}_e, \boldsymbol{\omega}_e)$. The perturbed vector fields $\hat{\boldsymbol{f}}_1$ and $\hat{\boldsymbol{f}}_2$ take the following form:

$$\hat{\boldsymbol{f}}_1(\boldsymbol{x}) = 0.5[\boldsymbol{E}(\boldsymbol{q}_e) - h\boldsymbol{I}_3]\boldsymbol{\omega}_e,$$

$$\hat{\boldsymbol{f}}_2(\boldsymbol{x}, t) = \boldsymbol{J}^{-1}\Xi(\boldsymbol{\omega}_e, \bar{\boldsymbol{\omega}}_d(t))\boldsymbol{\omega}_e - k_1\boldsymbol{J}^{-1}\bar{\boldsymbol{\kappa}}(h\boldsymbol{Q}_e, 1-\alpha_1).$$

Equations (20) and (21) take the same form as (2). Choose a Lyapunov candidate function as

$$\bar{V}_1(\boldsymbol{x}) = \frac{1}{2}\boldsymbol{\omega}_e^T \boldsymbol{J}\boldsymbol{\omega}_e + \frac{2k_1}{1+\alpha_1}\|\boldsymbol{q}_e\|^{1+\alpha_1}.$$

Note that $\mathrm{d}(\|\boldsymbol{q}_e\|^{1+\alpha_1})/\mathrm{d}t = (1+\alpha_1)\boldsymbol{\kappa}_0^T(\tilde{\boldsymbol{Q}}, 1-\alpha_1)\dot{\boldsymbol{q}}_e$. It can be proven by analysis similar to Steps 1) and 2) that the reduced system obtained from (20) and (21) by removing $\hat{\boldsymbol{f}}_1$ and $\hat{\boldsymbol{f}}_2$ is asymptotically stable with respect to the equilibria in $E_1$. Construct a dilation $\Delta_\varepsilon^r$ such that $\Delta_\varepsilon^r \boldsymbol{x} = (\varepsilon^{r_1}\boldsymbol{q}_e, \varepsilon^{r_2}\boldsymbol{\omega}_e)$, where $r_1 = -2k/(1-\alpha_2)$ and $r_2 = -(1+\alpha_2)k/(1-\alpha_2)$ for any $k < 0$. Recalling

$0 < \alpha_1 < 1$, $\alpha_2 = 2\alpha_1/(1+\alpha_1)$, it can be verified that the reduced system is homogeneous of degree $k$ with respect to $\Delta_\varepsilon^r$. In addition, recognizing $\boldsymbol{E}(\varepsilon^{r_1}\boldsymbol{q}_e) - h\boldsymbol{I}_3 \approx \mathcal{O}(\varepsilon^{r_1})\boldsymbol{q}_e^\times$ [18] and Property 2, it can be computed that $\hat{\boldsymbol{f}}_1(\Delta_\varepsilon^r \boldsymbol{x})/\varepsilon^{r_1+k} \to 0$ as $\varepsilon \to 0$. In addition, $\hat{\boldsymbol{f}}_2(\Delta_\varepsilon^r \boldsymbol{x}, t)/\varepsilon^{r_2+k} \to 0$ holds uniformly for all $t \geq 0$ as $\varepsilon \to 0$ since $\bar{\boldsymbol{\omega}}_d(t)$ is uniformly bounded according to Assumption 2.1. Lemma 2.1 can then be used to confirm uniform local finite-time stability of the equilibria in $E_1$.

Summarizing the analysis above, it is seen that the compact set $E_1$ is UGFTS. ■

### 3.2 GFTAC with Biased Angular Velocity Measurements

In this section it is assumed that the angular velocity measurements are corrupted by an unknown constant bias $\boldsymbol{b} \in \mathbb{R}^3$, i.e., the measured angular velocity is given by $\boldsymbol{\omega}_m = \boldsymbol{\omega} + \boldsymbol{b}$. In order to recover the true velocity information, an observer is designed to provide an estimate for the bias, denoted by $\hat{\boldsymbol{b}} \in \mathbb{R}^3$. The bias estimation error is defined as $\tilde{\boldsymbol{b}} = \boldsymbol{b} - \hat{\boldsymbol{b}}$. Introduce a virtual frame $\mathcal{F}_\mathcal{E}$ as an estimate of $\mathcal{F}_\mathcal{B}$ and its attitude relative to $\mathcal{F}_\mathcal{I}$ is denoted by $\boldsymbol{Q}_{EI} \in \mathbb{S}^3$. Let $\tilde{\boldsymbol{Q}} = [\tilde{q}_0, \tilde{\boldsymbol{q}}^T]^T \in \mathbb{S}^3$ represent the attitude error of $\mathcal{F}_\mathcal{B}$ relative to $\mathcal{F}_\mathcal{E}$ and by quaternion multiplication it follows that $\tilde{\boldsymbol{Q}} = \boldsymbol{Q}_{EI}^* \otimes \boldsymbol{Q}$.

Similarly to (13) and (14), define the following sets

$$C_2 = \{\boldsymbol{x}_2 \in M : \tilde{h}\tilde{q}_0 \geq -\delta\}, \tag{22}$$

$$D_2 = \{\boldsymbol{x}_2 \in M : \tilde{h}\tilde{q}_0 \leq -\delta\}, \tag{23}$$

where $\boldsymbol{x}_2 = (\tilde{\boldsymbol{Q}}, \tilde{\boldsymbol{b}}, \tilde{h}) \in M$ and $\tilde{h} \in \mathbb{H}$ is a switching variable associated with the hybrid observer below. The observer dynamics on the flow set $C_2$ is then designed as

$$\dot{\boldsymbol{Q}}_{EI} = \frac{1}{2}\boldsymbol{Q}_{EI} \otimes [\boldsymbol{R}^T(\tilde{\boldsymbol{Q}})(\boldsymbol{\omega}_m - \hat{\boldsymbol{b}} + \mu_1 \boldsymbol{\kappa}_1(\tilde{h}\tilde{\boldsymbol{Q}}, 1-\beta_1))], \tag{24}$$

$$\dot{\hat{\boldsymbol{b}}} = -\mu_2 \boldsymbol{\kappa}_1(\tilde{h}\tilde{\boldsymbol{Q}}, 1-\beta_2), \tag{25}$$

where $\mu_1, \mu_2 > 0$, $0.5 < \beta_1 < 1$, and $\beta_2 = 2\beta_1 - 1$. Note that $\dot{\tilde{h}} = 0$ during flows. When $\boldsymbol{x}_2 \in D_2$, the states follow the switching law $(\boldsymbol{Q}_{EI}^+, \hat{\boldsymbol{b}}^+, \tilde{h}^+) = (\boldsymbol{Q}_{EI}, \hat{\boldsymbol{b}}, \overline{\mathrm{sgn}}(\tilde{q}_0))$.

With the above observer, the equations of the estimation error are obtained as

$$\begin{cases} \dot{\tilde{\boldsymbol{Q}}} = \frac{1}{2}\tilde{\boldsymbol{Q}} \otimes [-\tilde{\boldsymbol{b}} - \mu_1 \boldsymbol{\kappa}_1(\tilde{h}\tilde{\boldsymbol{Q}}, 1-\beta_1))] \\ \dot{\tilde{\boldsymbol{b}}} = \mu_2 \boldsymbol{\kappa}_1(\tilde{h}\tilde{\boldsymbol{Q}}, 1-\beta_2) \end{cases}, \ \boldsymbol{x}_2 \in C_2, \tag{26}$$

$$\boldsymbol{x}_2^+ = G_2(\boldsymbol{x}_2) = (\tilde{\boldsymbol{Q}}, \tilde{\boldsymbol{b}}, \overline{\mathrm{sgn}}(\tilde{q}_0)), \ \boldsymbol{x}_2 \in D_2. \tag{27}$$

The closed-loop behavior of the observer is now described in the following theorem.

**Theorem 3.2**. Consider the hybrid control system given by (22), (23), (26) and (27) with $\mu_1, \mu_2 > 0$, $0.5 < \beta_1 < 1$, $\beta_2 = 2\beta_1 - 1$, and $\delta \in (0,1)$. Then, the compact set $E_2 = \{\boldsymbol{x}_2 \in M : \tilde{\boldsymbol{Q}} = \tilde{h}\mathbf{1}, \tilde{\boldsymbol{b}} = 0\}$ is GFTS.

**Proof**. This theorem can be proven in a manner similar to that of Theorem 3.1. Consider the Lyapunov candidate function

$$V_2(\boldsymbol{x}_2) = \frac{1}{2}\tilde{\boldsymbol{b}}^T\tilde{\boldsymbol{b}} + \frac{2\mu_2}{1+\beta_1}\varphi(\tilde{h}\tilde{q}_0, 1+\beta_1).$$

It can be readily computed that $\dot{V}_2(\boldsymbol{x}_2) = -\mu_1\mu_2 \left\| \boldsymbol{\kappa}_1(\tilde{h}\tilde{\boldsymbol{Q}}, 1-\beta_1) \right\|^2 \leq 0$ for $\boldsymbol{x}_2 \in C_2$ and $V_2(G_2(\boldsymbol{x}_2)) - V_2(\boldsymbol{x}_2) \leq -\sigma_2 < 0$ for $\boldsymbol{x}_2 \in D_2$, where $\sigma_2 = -2\mu_2\rho(-\delta, 1+\beta_1)/(1+\beta_1) > 0$. Hence, $E_2$ is Lyapunov stable and the number of jumps is uniformly bounded. In addition, it can be shown in a spirit analogous to the analysis in Step 2) of the proof of Theorem 3.1 that $\tilde{\boldsymbol{Q}}(t)$, $\tilde{\boldsymbol{b}}(t)$, $\dot{V}_2(t)$, and $\boldsymbol{\kappa}_1(\tilde{h}\tilde{\boldsymbol{Q}}(t), 1-\beta_1)$ are all uniformly continuous. It follows by invoking Barbalat's lemma that $\tilde{\boldsymbol{Q}}(t) \to \tilde{h}\mathbf{1}$ and $\tilde{\boldsymbol{b}}(t) \to 0$, and thus $E_2$ is GAS.

Next, (26) is decomposed into

$$\dot{\tilde{\boldsymbol{q}}} = -0.5\tilde{h}\tilde{\boldsymbol{b}} - 0.5\mu_1 \boldsymbol{\kappa}_0(\tilde{\boldsymbol{Q}}, 1-\beta_1) + \hat{\boldsymbol{f}}_1(\boldsymbol{x}), \tag{28}$$

$$\dot{\tilde{\boldsymbol{b}}} = \mu_2 \boldsymbol{\kappa}_0(\tilde{h}\tilde{\boldsymbol{Q}}, 1-\beta_2) + \hat{\boldsymbol{f}}_2(\boldsymbol{x}), \tag{29}$$

where $\boldsymbol{x} = (\tilde{\boldsymbol{q}}, \tilde{\boldsymbol{b}})$ and

$$\hat{\boldsymbol{f}}_1(\boldsymbol{x}) = -0.5[\boldsymbol{E}(\tilde{\boldsymbol{q}}) - \tilde{h}\boldsymbol{I}_3]\tilde{\boldsymbol{b}} - 0.5\mu_1[(\tilde{q}_0 - \tilde{h})\boldsymbol{\kappa}_1(\tilde{h}\tilde{\boldsymbol{Q}}, 1-\beta_1) + \tilde{h}\bar{\boldsymbol{\kappa}}(\tilde{h}\tilde{\boldsymbol{Q}}, 1-\beta_1)], \tag{30}$$

$$\hat{\boldsymbol{f}}_2(\boldsymbol{x}) = \mu_2 \bar{\boldsymbol{\kappa}}(\tilde{h}\tilde{\boldsymbol{Q}}, 1-\beta_2)). \tag{31}$$

Note that $\tilde{\boldsymbol{q}}^{\times}\boldsymbol{\kappa}_1(\tilde{h}\tilde{\boldsymbol{Q}},\cdot)=0$ and $\boldsymbol{\kappa}_0(\tilde{\boldsymbol{Q}},\cdot)=\tilde{h}\boldsymbol{\kappa}_0(\tilde{h}\tilde{\boldsymbol{Q}},\cdot)$ are utilized to obtain (30). Discarding $\hat{\boldsymbol{f}}_1(\boldsymbol{x})$ and $\hat{\boldsymbol{f}}_2(\boldsymbol{x})$ from (28) and (29) yields a reduced system, which can be proven to be asymptotically stable by examining the Lyapunov function candidate $\bar{V}_2(\boldsymbol{x})=\tilde{\boldsymbol{b}}^T\tilde{\boldsymbol{b}}/2+2\mu_2\|\tilde{\boldsymbol{q}}\|^{1+\beta_1}/(1+\beta_1)$. In addition, the reduced system is homogeneous of degree $k<0$ with respect to a dilation $\Delta_\varepsilon^{\boldsymbol{r}}\boldsymbol{x}=(\varepsilon^{r_1}\tilde{\boldsymbol{q}},\varepsilon^{r_2}\tilde{\boldsymbol{b}})$, where $r_1=-k/(1-\beta_1)$ and $r_2=-\beta_1 k/(1-\beta_1)$. Using $\boldsymbol{E}(\varepsilon^{r_1}\tilde{\boldsymbol{q}})-\tilde{h}\boldsymbol{I}_3\approx\mathcal{O}(\varepsilon^{r_1})\tilde{\boldsymbol{q}}^{\times}$, $\tilde{q}_0(\varepsilon^{r_1}\tilde{\boldsymbol{q}})-\tilde{h}\approx\mathcal{O}(\varepsilon^{2r_1})\|\tilde{\boldsymbol{q}}\|^2$ [18] and Property 2, further computations show that $\hat{\boldsymbol{f}}_i(\Delta_\varepsilon^{\boldsymbol{r}}\boldsymbol{x})/\varepsilon^{r_i+k}\to 0$ as $\varepsilon\to 0$, $i\in\mathbb{I}_2$. Invoking Lemma 2.1 ensures local finite-time stability of $E_2$, which when being synthesized with the previous GAS result leads to GFTS of $E_2$. ■

**Remark 3.3**. By setting $\beta_1=\beta_2=1$ and $\tilde{h}=\overline{\operatorname{sgn}}(\tilde{q}_0(t))$ in (24) and (25), the discontinuous exponential observer derived in [5] can be retrieved and the proposed hybrid observer not only avoids the chattering due to measurement noise but also provides the true bias in finite time. More importantly, this latter property enables one to obtain a GFTAC by combining the hybrid observer and the previous full-state controller. A detailed analysis is presented in the sequel.

With a bias estimate $\hat{\boldsymbol{b}}$, one can define $\hat{\boldsymbol{\omega}}=\boldsymbol{\omega}_m-\hat{\boldsymbol{b}}$ and $\hat{\boldsymbol{\omega}}_e=\hat{\boldsymbol{\omega}}-\bar{\boldsymbol{\omega}}_d$ as the estimates of true velocity and velocity tracking error, respectively. Letting $\hat{\boldsymbol{x}}_1=(\boldsymbol{Q}_e,\hat{\boldsymbol{\omega}}_e,h)$, a certainty-equivalence controller $\boldsymbol{u}(\hat{\boldsymbol{x}}_1)$ can be obtained from $\boldsymbol{u}(\boldsymbol{x}_1)$ defined in (13)-(16) by substituting $\hat{\boldsymbol{\omega}}_e$ for $\boldsymbol{\omega}_e$. Since the proposed observer ensures finite-time recovery of the true bias, $\boldsymbol{u}(\hat{\boldsymbol{x}}_1)$ is restored to $\boldsymbol{u}(\boldsymbol{x}_1)$ within a finite time if $(\boldsymbol{Q}_e,\boldsymbol{\omega}_e)$ does not escape in finite time. Consequently, $\boldsymbol{u}(\hat{\boldsymbol{x}}_1)$ can also lead to finite-time convergence of $(\boldsymbol{Q}_e,\boldsymbol{\omega}_e)$ to $(h\mathbf{1},0)$, according to Theorem 3.1. The settling time of $\boldsymbol{u}(\hat{\boldsymbol{x}}_1)$ is bounded by the addition of the settling times of the preceding bias observer and the finite-time controller $\boldsymbol{u}(\boldsymbol{x}_1)$, both of them relying on initial values of the system states.

To demonstrate that there is no finite escape time for the closed-loop system, consider again $V_1(\boldsymbol{x}_1)$ given in (17). Its time derivative under the effect of $\boldsymbol{u}(\hat{\boldsymbol{x}}_1)$ now becomes

$$\begin{aligned}\dot{V}_1&=-k_2\boldsymbol{\omega}_e^T\operatorname{sat}_{\alpha_2}(\boldsymbol{\omega}_e)+k_2\boldsymbol{\omega}_e^T[\operatorname{sat}_{\alpha_2}(\boldsymbol{\omega}_e)-\operatorname{sat}_{\alpha_2}(\hat{\boldsymbol{\omega}}_e)]\\&\le -k_2\boldsymbol{\omega}_e^T\operatorname{sat}_{\alpha_2}(\boldsymbol{\omega}_e)+k_2\sum_{i=1}^{3}|\omega_{ei}|\left|\operatorname{sat}_{\alpha_2}(\omega_{ei})-\operatorname{sat}_{\alpha_2}(\hat{\omega}_{ei})\right|\\&\le -k_2\boldsymbol{\omega}_e^T\operatorname{sat}_{\alpha_2}(\boldsymbol{\omega}_e)+k_2\sum_{i=1}^{3}2|\omega_{ei}|\le 3k_2+k_2\|\boldsymbol{\omega}_e\|^2\\&\le 3k_2+c_0V_1\end{aligned}\quad,\ \boldsymbol{x}_1\in C_1, \tag{32}$$

where $c_0=2k_2/\lambda_{\mathrm{m}}(\boldsymbol{J})$ and $\lambda_{\mathrm{m}}(\cdot)$ represents the minimum eigenvalue of a square matrix. Equation (32) implies that $V_1(t)\le V_1(t_0)+3k_2(t-t_0)+\int_{t_0}^{t}c_0V_1(s)\,\mathrm{d}s$. Employing the Gronwall-Bellman inequality [26] leads to

$$V_1(t)\le c_1\exp(c_0(t-t_0))-3k_2/c_0,\ \boldsymbol{x}_1\in C_1,\ t\ge t_0, \tag{33}$$

where $c_1=V_1(t_0)+3k_2/c_0$. In addition, note that $(\boldsymbol{Q}_e,\boldsymbol{\omega}_e)$ remains continuous when jumps occur, i.e., $\boldsymbol{x}_1\in D_1$. Equation (33) then implies that $(\boldsymbol{Q}_e,\boldsymbol{\omega}_e)$ remains bounded (and thus cannot escape) in finite time under $\boldsymbol{u}(\hat{\boldsymbol{x}}_1)$. The so-called "separation principle" is now obtained. The preceding analysis leads to the following proposition:

**Proposition 3.1**. Consider the certainty-equivalence pair of the hybrid controller given by (13)-(16) with $\boldsymbol{\omega}_e$ replaced by $\hat{\boldsymbol{\omega}}_e=\hat{\boldsymbol{\omega}}-\bar{\boldsymbol{\omega}}_d=\boldsymbol{\omega}_m-\hat{\boldsymbol{b}}-\bar{\boldsymbol{\omega}}_d$, where $\hat{\boldsymbol{b}}$ is generated from the observer (22)-(25). Then, the attitude tracking error $(\boldsymbol{Q}_e,\boldsymbol{\omega}_e)$ is uniformly bounded and globally converges to $(h\mathbf{1},0)$ in finite time.

### 3.3 GFTAC with Attitude-Only Measurements

Next, only the quaternion attitude information is assumed to be available. Again, consider the virtual frame $\mathcal{F}_{\mathcal{E}}$ and denote by $\boldsymbol{Q}_{ED}\in\mathbb{S}^3$ its attitude relative to $\mathcal{F}_{\mathcal{D}}$. The attitude error of $\mathcal{F}_{\mathcal{B}}$ relative to $\mathcal{F}_{\mathcal{E}}$ is given by $\tilde{\boldsymbol{Q}}=\boldsymbol{Q}_{ED}^*\otimes\boldsymbol{Q}_e$. In the following, an auxiliary system is constructed to dynamically update $\boldsymbol{Q}_{ED}$ such that $\tilde{\boldsymbol{Q}}$ provides damping, which is otherwise provided by velocity feedback.

Denote $\tilde{M}\triangleq\mathbb{S}^3\times M\times\mathbb{H}$ and $\boldsymbol{x}_3=(\tilde{\boldsymbol{Q}},\boldsymbol{x}_1,\tilde{h})\in\tilde{M}$, where $\tilde{h}\in\mathbb{H}$ is the switching variable accompanying $\boldsymbol{Q}_{ED}$ (or $\tilde{\boldsymbol{Q}}$). The objective is to stabilize the set $E_3=\{\boldsymbol{x}_3\in\tilde{M}:\boldsymbol{x}_1\in E_1,\tilde{\boldsymbol{Q}}=\tilde{h}\mathbf{1}\}$ in finite time. To this end, define the flow and jump sets as follows:

$$C_3=\{\boldsymbol{x}_3\in\tilde{M}:hq_{e0}\ge-\delta\text{ and }\tilde{h}\tilde{q}_0\ge-\delta\}, \tag{34}$$

$$D_3=\{\boldsymbol{x}_3\in\tilde{M}:hq_{e0}\le-\delta\text{ or }\tilde{h}\tilde{q}_0\le-\delta\}. \tag{35}$$

The quaternion filter and the control torque are now designed as

$$\dot{\boldsymbol{Q}}_{ED}=\frac{1}{2}\boldsymbol{Q}_{ED}\otimes[k_3\boldsymbol{R}^T(\tilde{\boldsymbol{Q}})\boldsymbol{\kappa}_1(\tilde{h}\tilde{\boldsymbol{Q}},1-\alpha_3))],\ \boldsymbol{x}_3\in C_3, \tag{36}$$

$$\boldsymbol{u}(\boldsymbol{x}_3)=\boldsymbol{u}_d-k_1\boldsymbol{\kappa}_1(h\boldsymbol{Q}_e,1-\alpha_1)-k_2\boldsymbol{\kappa}_1(\tilde{h}\tilde{\boldsymbol{Q}},1-\alpha_1),\ \boldsymbol{x}_3\in C_3, \tag{37}$$

where $k_i > 0$, $i \in \mathbb{I}_3$, $0.5 < \alpha_3 < 1$, and $\alpha_1 = 2\alpha_3 - 1$. Note that $\dot{h} = \dot{\tilde{h}} = 0$ during flows. If $\boldsymbol{x}_3 \in D_3$, the states jump to

$$\boldsymbol{x}_3^+ = G_3(\boldsymbol{x}_3) = (\tilde{\boldsymbol{Q}}, \boldsymbol{Q}_e, \boldsymbol{\omega}_e, \overline{\mathrm{sgn}}(q_{e0}), \overline{\mathrm{sgn}}(\tilde{q}_0)),\ \boldsymbol{x}_3 \in D_3. \tag{38}$$

The time derivative of $\tilde{\boldsymbol{Q}}$ can be obtained by employing (36):

$$\dot{\tilde{\boldsymbol{Q}}} = \frac{1}{2}\tilde{\boldsymbol{Q}} \otimes [\boldsymbol{\omega}_e - k_3 \boldsymbol{\kappa}_1(\tilde{h}\tilde{\boldsymbol{Q}}, 1-\alpha_3))],\ \boldsymbol{x}_3 \in C_3. \tag{39}$$

The following theorem states the stability properties of the overall closed-loop system.

**Theorem 3.3**. Consider the hybrid control system given by (10), (11), (34), (35), and (37)-(39) with $k_i > 0$, $i \in \mathbb{I}_3$, $0.5 < \alpha_3 < 1$, $\alpha_1 = 2\alpha_3 - 1$, and $\delta \in (0,1)$. Then, the compact set $E_3$ is UGFTS.

**Remark 3.4**. Invoking Assumption 2.1 and Property 1, it follows that $\boldsymbol{u}(\boldsymbol{x}_3)$ is bounded by $|u_i(\boldsymbol{x}_3)| < k_1 + k_2 + (\varpi_1^2 + \varpi_2)\|\boldsymbol{J}\|$, $i \in \mathbb{I}_3$. In addition, the proposed output-feedback control law reduces to the hybrid output-feedback controller in [7] by setting $\alpha_1 = \alpha_2 = 1$ and the smooth output-feedback controller in [3] by further setting $h \equiv 1$ and replacing $k_3$ with a positive-definite matrix. These two methods both achieve asymptotic convergence, rather than the more desirable finite-time convergence as obtained in this paper.

**Remark 3.5**. All the above hybrid control schemes overcome the topological constraint on the attitude manifold, thus ensuring global attitude tracking, by means of discontinuous dynamics when system states reside on the jump sets. Such discontinuities are not a major problem for the proposed bias observer and quaternion filter because they are numerically implemented. The resultant discontinuous command torque is compatible with physical actuators in the on/off mode of operation such as thrusters, but cannot be implemented by actuators that can only provide continuous inputs.

**Proof of Theorem 3.3**. The proof is similar to that of Theorem 3.1 with crucial modifications being sketched below. Consider the Lyapunov candidate function

$$V_3(\boldsymbol{x}_3) = V_1(\boldsymbol{x}_1) + \frac{2k_2}{1+\alpha_3}\varphi(\tilde{h}\tilde{q}_0, 1+\alpha_3).$$

Direct computation shows that $\dot{V}_3(\boldsymbol{x}_3) = -k_1 k_2 k_3 \left\|\boldsymbol{\kappa}_1(\tilde{h}\tilde{\boldsymbol{Q}}, 1-\alpha_3)\right\|^2 \le 0$ for $\boldsymbol{x}_3 \in C_3$ and

$$V_3(G_3(\boldsymbol{x}_3)) - V_3(\boldsymbol{x}_3) = \frac{2k_1}{1+\alpha_1}\rho(hq_{e0}, 1+\alpha_1) + \frac{2k_2}{1+\alpha_3}\rho(\tilde{h}\tilde{q}_0, 1+\alpha_3),$$

for $\boldsymbol{x}_3 \in D_3$, where it follows that $\rho(hq_{e0}, \cdot) \le 0$ and $\rho(\tilde{h}\tilde{q}_0, \cdot) \le 0$, and at least one of them is negative. Hence, $V_3(G_3(\boldsymbol{x}_3)) - V_3(\boldsymbol{x}_3) < 0$ and $E_3$ is uniformly Lyapunov stable and the number of jumps is uniformly bounded. In addition, it can be shown in a spirit similar to the analysis in Step 2) of the proof of Theorem 3.1 that $\boldsymbol{Q}_e(t)$, $\tilde{\boldsymbol{Q}}(t)$, $\boldsymbol{\omega}_e(t)$, $\dot{V}_3(t)$, $\dot{\boldsymbol{\omega}}_e(t)$ and $\boldsymbol{\kappa}_1(\tilde{h}\tilde{\boldsymbol{Q}}(t), 1-\beta_1)$ are all uniformly continuous. It follows by invoking Barbalat's lemma that $\tilde{\boldsymbol{Q}}(t) \to \tilde{h}\boldsymbol{1}$, $\boldsymbol{Q}_e(t) \to h\boldsymbol{1}$ and $\boldsymbol{\omega}_e(t) \to 0$, and thus $E_3$ is GAS.

Similarly to the proof Theorems 3.1 and 3.2, the closed-loop flow dynamics can also be written as

$$\dot{\tilde{\boldsymbol{q}}} = 0.5\tilde{h}\boldsymbol{\omega}_e - 0.5k_3\boldsymbol{\kappa}_0(\tilde{\boldsymbol{Q}}, 1-\alpha_3) + \hat{\boldsymbol{f}}_1(\boldsymbol{x}), \tag{40}$$

$$\dot{\boldsymbol{q}}_e = 0.5h\boldsymbol{\omega}_e + \hat{\boldsymbol{f}}_2(\boldsymbol{x}), \tag{41}$$

$$\dot{\boldsymbol{\omega}}_e = -\boldsymbol{J}^{-1}[k_1\boldsymbol{\kappa}_0(h\boldsymbol{Q}_e, 1-\alpha_1) + k_2\boldsymbol{\kappa}_0(\tilde{h}\tilde{\boldsymbol{Q}}, 1-\alpha_1)] + \hat{\boldsymbol{f}}_3(\boldsymbol{x}, t), \tag{42}$$

where $\boldsymbol{x} = (\tilde{\boldsymbol{q}}, \boldsymbol{q}_e, \boldsymbol{\omega}_e)$ and

$$\hat{\boldsymbol{f}}_1(\boldsymbol{x}) = 0.5[\boldsymbol{E}(\tilde{\boldsymbol{q}}) - \tilde{h}\boldsymbol{I}_3]\boldsymbol{\omega}_e - 0.5k_3(\tilde{q}_0 - \tilde{h})\boldsymbol{\kappa}_1(\tilde{h}\tilde{\boldsymbol{Q}}, 1-\alpha_3) - 0.5k_3\tilde{h}\bar{\boldsymbol{\kappa}}(\tilde{h}\tilde{\boldsymbol{Q}}, 1-\alpha_3),$$

$$\hat{\boldsymbol{f}}_2(\boldsymbol{x}) = 0.5[\boldsymbol{E}(\boldsymbol{q}_e) - h\boldsymbol{I}_3]\boldsymbol{\omega}_e,$$

$$\hat{\boldsymbol{f}}_3(\boldsymbol{x}, t) = \boldsymbol{J}^{-1}\Xi(\boldsymbol{\omega}_e, \bar{\boldsymbol{\omega}}_d(t))\boldsymbol{\omega}_e - \boldsymbol{J}^{-1}[k_1\bar{\boldsymbol{\kappa}}(h\boldsymbol{Q}_e, 1-\alpha_1) + k_2\bar{\boldsymbol{\kappa}}(\tilde{h}\tilde{\boldsymbol{Q}}, 1-\alpha_1)].$$

Omitting $\hat{\boldsymbol{f}}_i$, $i \in \mathbb{I}_3$, from (40)-(42) yields a reduced system. Construct a Lyapunov function candidate $\bar{V}_3(\boldsymbol{x}) = \bar{V}_1(\boldsymbol{x}) + 2k_3\|\tilde{\boldsymbol{q}}\|^{1+\alpha_3}/(1+\alpha_3)$ and dilation $\Delta_\varepsilon^r \boldsymbol{x} = (\varepsilon^{r_1}\tilde{\boldsymbol{q}}, \varepsilon^{r_2}\boldsymbol{q}_e, \varepsilon^{r_3}\boldsymbol{\omega}_e)$, where $r_2 = r_3 = -k/(1-\alpha_3)$ and $r_1 = -\alpha_3 k/(1-\alpha_3)$. It can then be shown that the reduced system is asymptotically stable and homogeneous of degree $k < 0$. In addition, $\hat{\boldsymbol{f}}_i(\Delta_\varepsilon^r \boldsymbol{x})/\varepsilon^{r_i+k} \to 0$, $i \in \mathbb{I}_2$, and $\hat{\boldsymbol{f}}_3(\Delta_\varepsilon^r \boldsymbol{x}, t)/\varepsilon^{r_3+k} \to 0$ hold uniformly for all $t \ge 0$ as $\varepsilon \to 0$. Finally, Lemma 2.1 is used to conclude the uniform finite-time stability of $E_3$ at the local and then global levels. ■

# 4 Simulations

In this section three numerical examples are presented to illustrate the performance of the proposed methods. The rigid body, with an inertia matrix of $\boldsymbol{J}=\mathrm{diag}\{15,20,10\}$ kg·m², is required to track a trajectory given by $\boldsymbol{Q}_d(0)=\mathbf{1}$ and $\boldsymbol{\omega}_d(t)=0.01\times[\sin(\omega_0 t),\sin(\omega_0 t),\sin(\omega_0 t)]^T$ rad/s, where $\omega_0=0.01$ rad/s. It is assumed that a disturbance of the form $\boldsymbol{d}(t)=2[\cos(0.1t),\cos(0.1t),-\sin(0.1t)]^T\times 10^{-2}$ N·m acts on the rigid body. The attitude measurement noise is modeled such that the eigenaxis $\boldsymbol{\eta}$ associated with the measured $\boldsymbol{q}$ is uniformly distributed within a spherical cone centered at the true eigenaxis with a half cone angle of 0.01 deg. In addition, the angular velocity measurement model is given by $\boldsymbol{\omega}_m(t)=\boldsymbol{\omega}(t)+\boldsymbol{b}(t)+\boldsymbol{v}(t)$ and $\dot{\boldsymbol{b}}(t)=\boldsymbol{v}_b(t)$, where the gyro noise $\boldsymbol{v}(t)\in\mathbb{R}^3$ and the bias noise $\boldsymbol{v}_b(t)\in\mathbb{R}^3$ are random, zero mean, Gaussian noise processes with standard deviations of 0.01 deg/s and 0.01 deg/s², respectively. Note that in the full-state measurement case (Example 1) only the gyro noise is considered while $\boldsymbol{b}(t)$ is set to zero. The initial attitude and angular velocity of the rigid body are set as $q_0(0)=0$, $\boldsymbol{q}(0)=[0.6,-0.8,0]^T$, and $\boldsymbol{\omega}(0)=[0.3,-0.4,0]^T$ rad/s. The maximum control torque is 5 N·m.

## 4.1 Example 1: Full-State Measurements

The gains for the GFTAC with full-state feedback are chosen as $k_1=1.1$, $k_2=4$, $\delta=0.3$, and $h(0)=1$. The cases of $\alpha_1=0.6,0.8,1$ ($\alpha_2=2\alpha_1/(1+\alpha_1)$) are examined. Note that $\alpha_1=1$ ($\alpha_2=1$) corresponds to the asymptotic controller in [7].

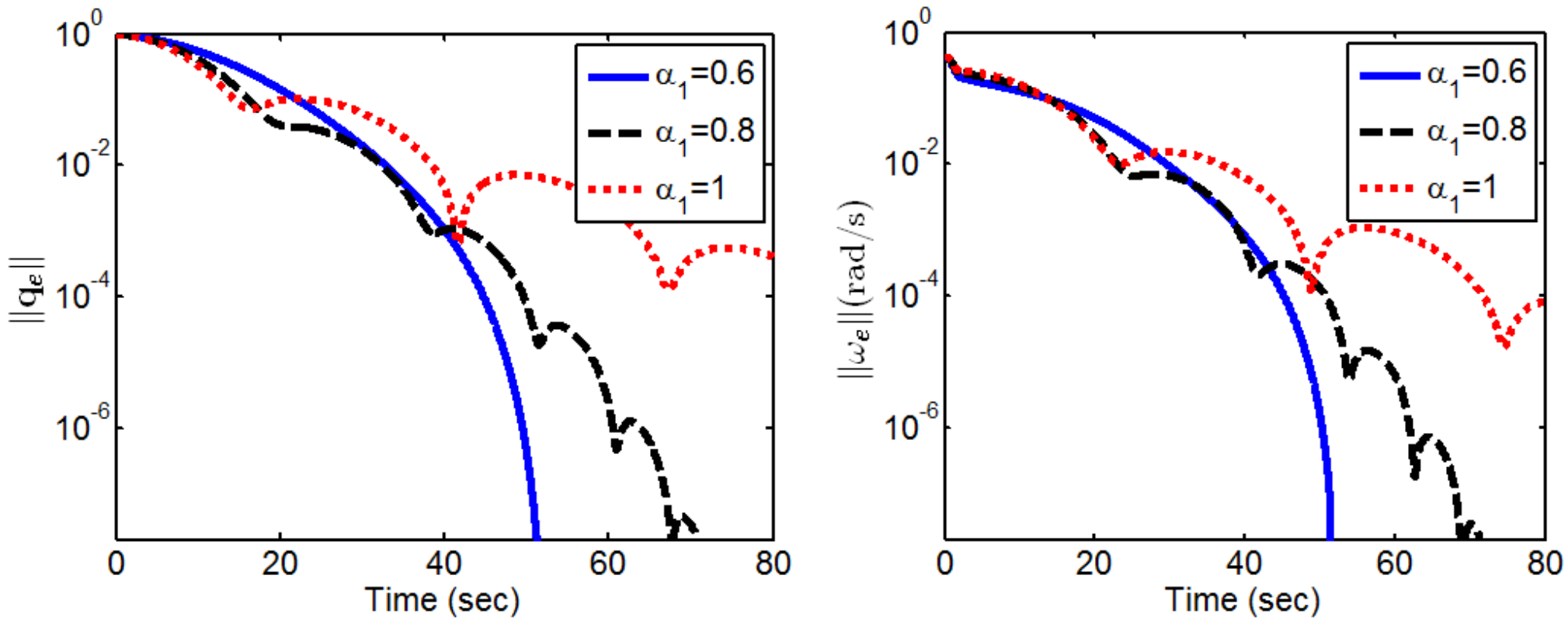

a) norm of $\boldsymbol{q}_e$ b) norm of $\boldsymbol{\omega}_e$

**Fig. 1 Time histories of a) $\|\boldsymbol{q}_e\|$ and b) $\|\omega_e\|$ with no uncertainties.**

Figure 1 plots the norms of $\boldsymbol{q}_e$ and $\boldsymbol{\omega}_e$ in the absence of any system uncertainty. In this case, finite-time convergence of the tracking error is achieved for $\alpha_1=0.6,0.8$ at about 55 s and 75 s respectively, which verifies the theoretical result in Theorem 3.1. In contrast, infinite convergence time is expected for the asymptotic method (i.e., $\alpha_1=1$). It can be seen that as the tracking error is close to zero, decreasing the value of $\alpha_1$ can significantly increase the convergence rate.

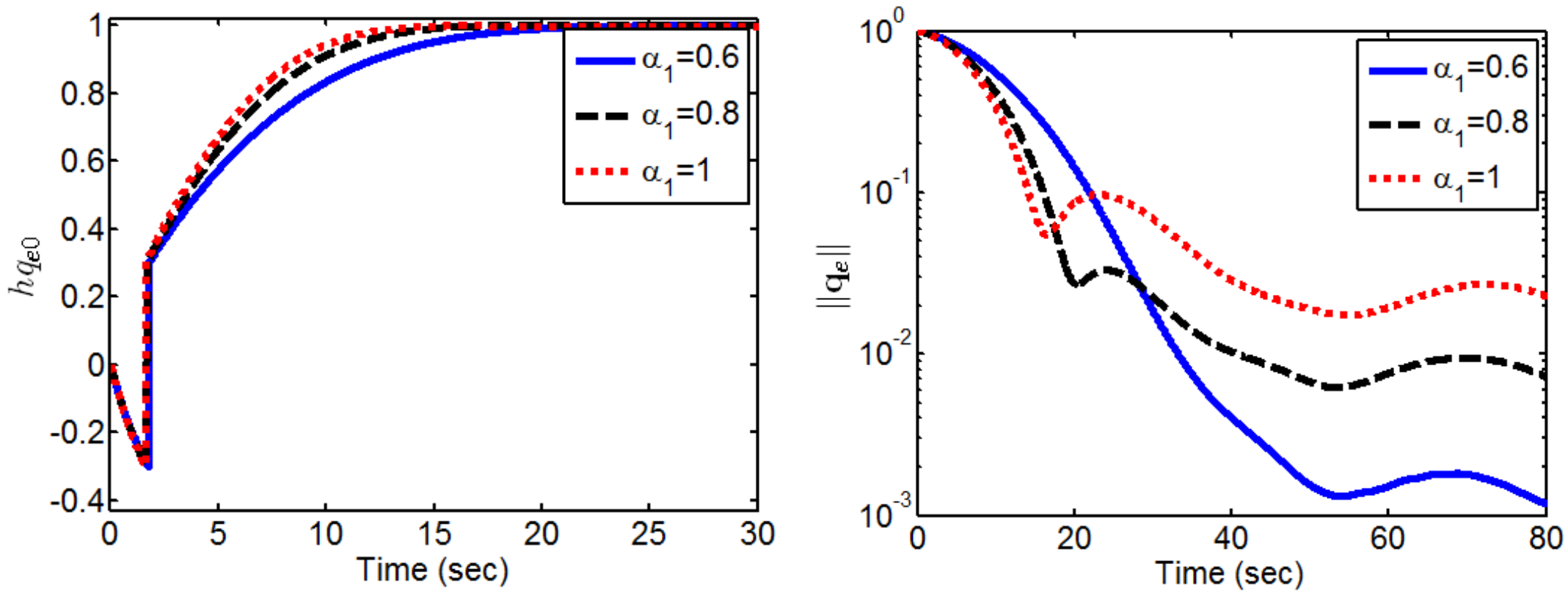

a) $hq_{e0}$ b) norm of $\boldsymbol{q}_e$

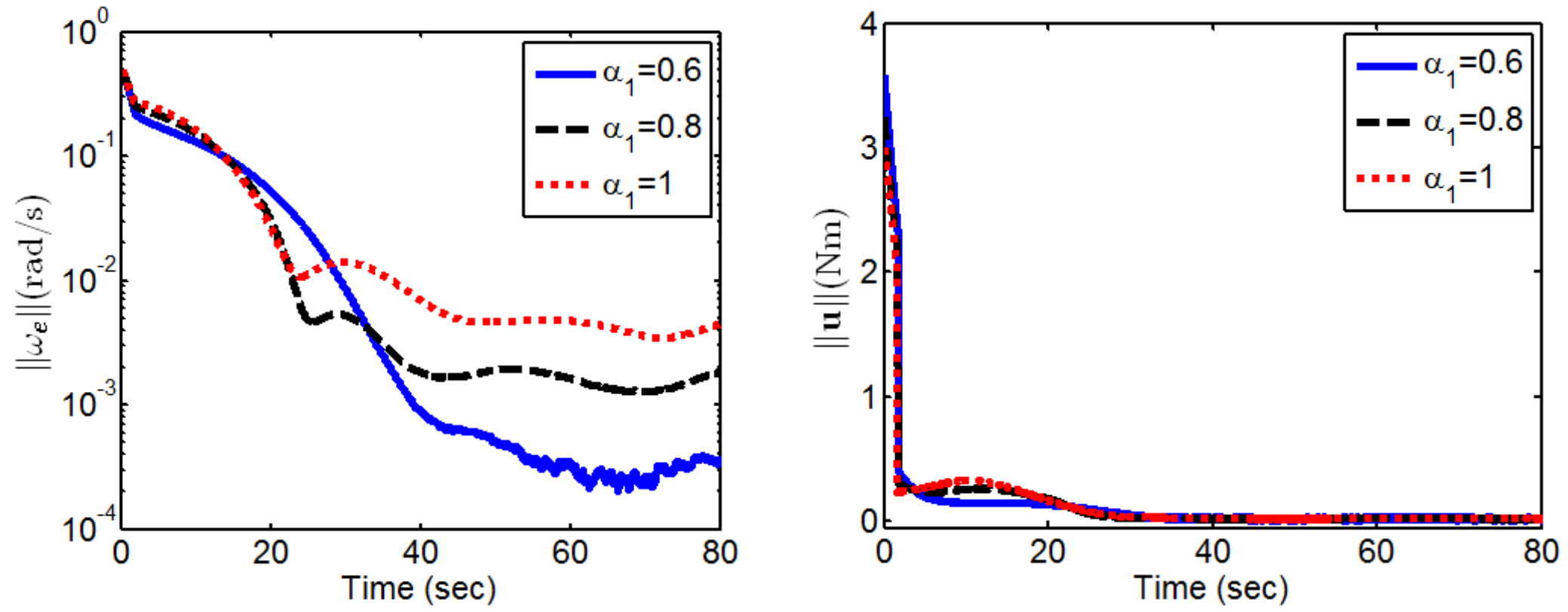

c) norm of $\boldsymbol{\omega}_e$ d) norm of $\boldsymbol{u}$

**Fig. 2 Simulation results for the GFTAC with full-state feedback and uncertainties.**

The simulation results with both external disturbances and measurement noise are given in Fig. 2. Figure 2a depicts the time histories of $hq_{e0}$. Although the control law commands $\boldsymbol{Q}_e$ toward $\mathbf{1}$ from the onset, the initial angular velocity drives $\boldsymbol{Q}_e$ toward $-\mathbf{1}$. As $\boldsymbol{Q}_e$ moves past the hysteresis gap $\delta = 0.3$ at about 1.5 s, $h$ jumps to $-1$ with robustness against measurement noise. For all selected values of $\alpha_1$, $(\boldsymbol{Q}_e, \boldsymbol{\omega}_e)$ is eventually stabilized to a small neighborhood of $(-\mathbf{1}, 0)$, as shown in Figs. 2b and 2c, and the control torques are within the saturation limit, as shown in Fig. 2d. The resultant steady-state accuracy for $\alpha_1 = 0.6, 0.8, 1$, however, is greatly different. It can be seen that the finite-time method achieves smaller tracking error and thus better robustness than the asymptotic method. In addition, decreasing the value of $\alpha_1$ further reduces the tracking error but increases the nonsmooth degree of control torques.

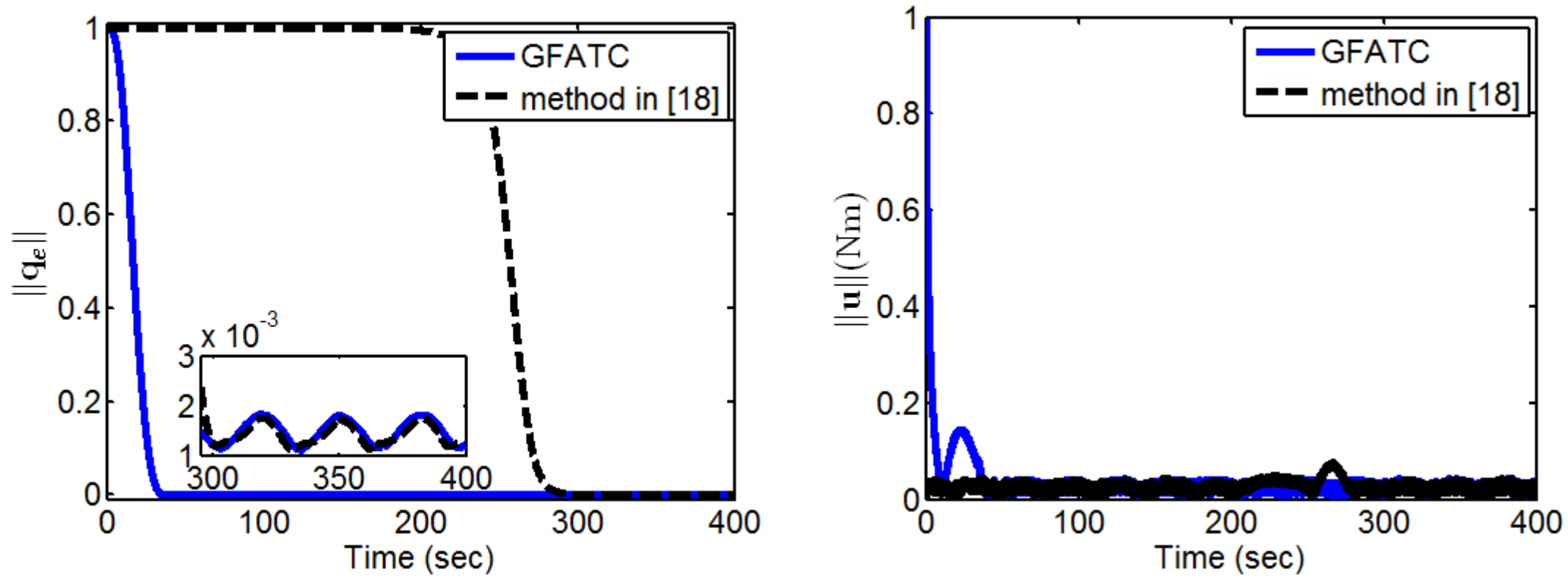

a) norm of $\boldsymbol{q}_e$ b) norm of $\boldsymbol{u}$

**Fig. 3 Performance of the GFTAC and the method in [18].**

Next, the GFTAC and the quaternion-based controller derived in [18] with AGFTS are simulated for a large angle maneuver, i.e., $\boldsymbol{Q}_d \equiv \mathbf{1}$ and $\boldsymbol{\omega}_d \equiv [0,0,0]^T$ rad/s. The disturbance and noise models remain the same while the initial conditions are renewed as $\boldsymbol{q}(0) = [1, 0, 0]^T$ and $\boldsymbol{\omega}(0) = [0,0,0]^T$ rad/s. In other words, the initial state resides on an undesired equilibrium of the method in [18]. We select $k_1 = 1.1$, $k_2 = 4$, and $\alpha_1 = 0.6$ for both the GFTAC and the method in [18], and $\delta = 0.3$ and $h(0) = 1$.

Figure 3 plots the time histories of $\|\boldsymbol{q}_e\|$ and $\|\boldsymbol{u}\|$ obtained by the two methods. We can see a drastic contrast in their convergence rates, although both methods ensure a similar steady-state tracking accuracy. The response of $\|\boldsymbol{q}_e\|$ by the method of [18] is rather sluggish during $[0, 240]$ s, because the applied control torque is close to zero around the initial state. As the attitude is gradually perturbed away from the undesired equilibrium by the external disturbance and measurement noise, the method of [18] takes effect and eventually maneuvers the rigid body to the target attitude. The total maneuver time is about 295 s. In contrast, the GFTAC achieves a faster maneuver with a convergence time of 35 s, attributable to its advantage of global stability.

### 4.2 Example 2: Biased Angular Velocity Measurements

In this example, assume that angular velocity measurements possess a bias of $\boldsymbol{b}(0) = [1, -5, 2]^T \times 10^{-2}$ rad/s. The gains and initial estimates for the proposed hybrid bias observer are selected as $\mu_1 = 0.33$, $\mu_2 = 0.12$, $\beta_1 = 0.75$, $\tilde{h}(0) = 1$, $\boldsymbol{Q}_{EI}(0) = \boldsymbol{Q}(0)$, and $\hat{\boldsymbol{b}}(0) = [0,0,0]^T$ rad/s while the gains for the certainty-equivalence controller $\boldsymbol{u}(\hat{\boldsymbol{x}}_1)$ remain the same as Example 1.

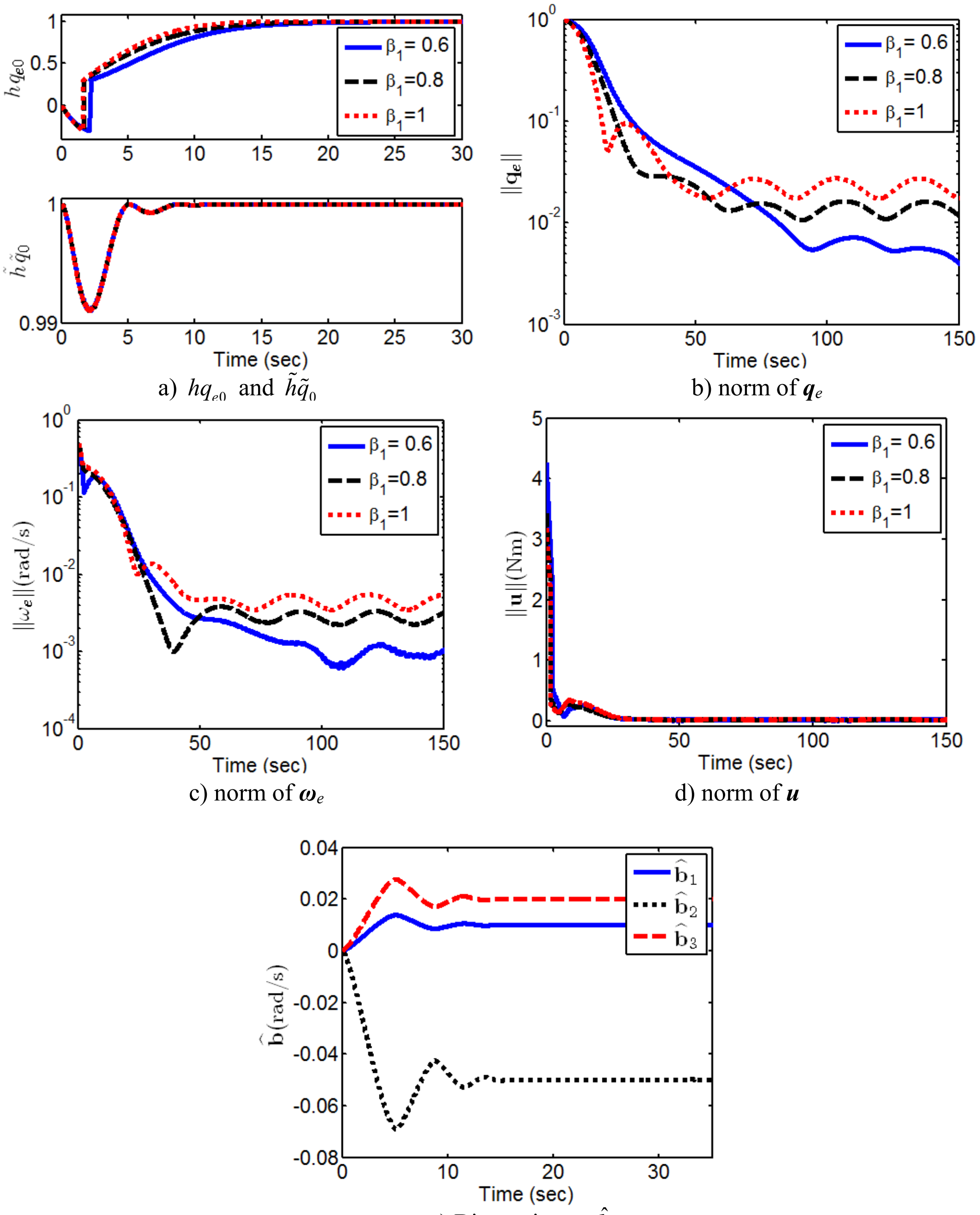

a) $hq_{e0}$ and $\tilde{h}\tilde{q}_0$ b) norm of $\boldsymbol{q}_e$

c) norm of $\boldsymbol{\omega}_e$ d) norm of $\boldsymbol{u}$

e) Bias estimate $\hat{\boldsymbol{b}}$

**Fig. 4 Simulation results for the GFTAC with biased angular velocity measurements and uncertainties.**

The simulation results are given in Fig. 4. As shown in Fig. 4a, $h$ switches once from 1 to $-1$ and chattering (i.e., multiple jumps occurring at the same time) is avoided since the hysteresis width $\delta = 0.3$ is larger than the noise bound. The sign of $\tilde{h}$ remains the same because $\tilde{q}_0$ is positive during the whole control phase. Despite the presence of external disturbances, measurement noise and gyro bias, the certainty-equivalence controller $\boldsymbol{u}(\hat{\boldsymbol{x}}_1)$ still stabilizes the tracking error to a small neighborhood of zero. Similarly to Example 1, smaller steady-state tracking error can be seen by decreasing the power gain $\alpha_1$ and control torques are below the saturation limit. For all selected values of $\alpha_1$, the bias observer exhibits the same behavior because it is independent of the controller $\boldsymbol{u}(\hat{\boldsymbol{x}}_1)$. As shown in Fig. 4e, the bias estimate $\hat{\boldsymbol{b}}$ converges approximately instead of exactly to the true bias due to the effect of measurement noise.

### 4.3 Example 3: Attitude-Only Measurements

In this part, only attitude measurements are assumed to be available. The control parameters for the velocity-free GFTAC are chosen as $k_1 = 1.2$, $k_2 = 2.4$, $k_3 = 1.1$, $\delta = 0.3$, $h(0) = \tilde{h}(0) = 1$, and $\boldsymbol{Q}_{ED}(0) = \boldsymbol{Q}_e(0)$.

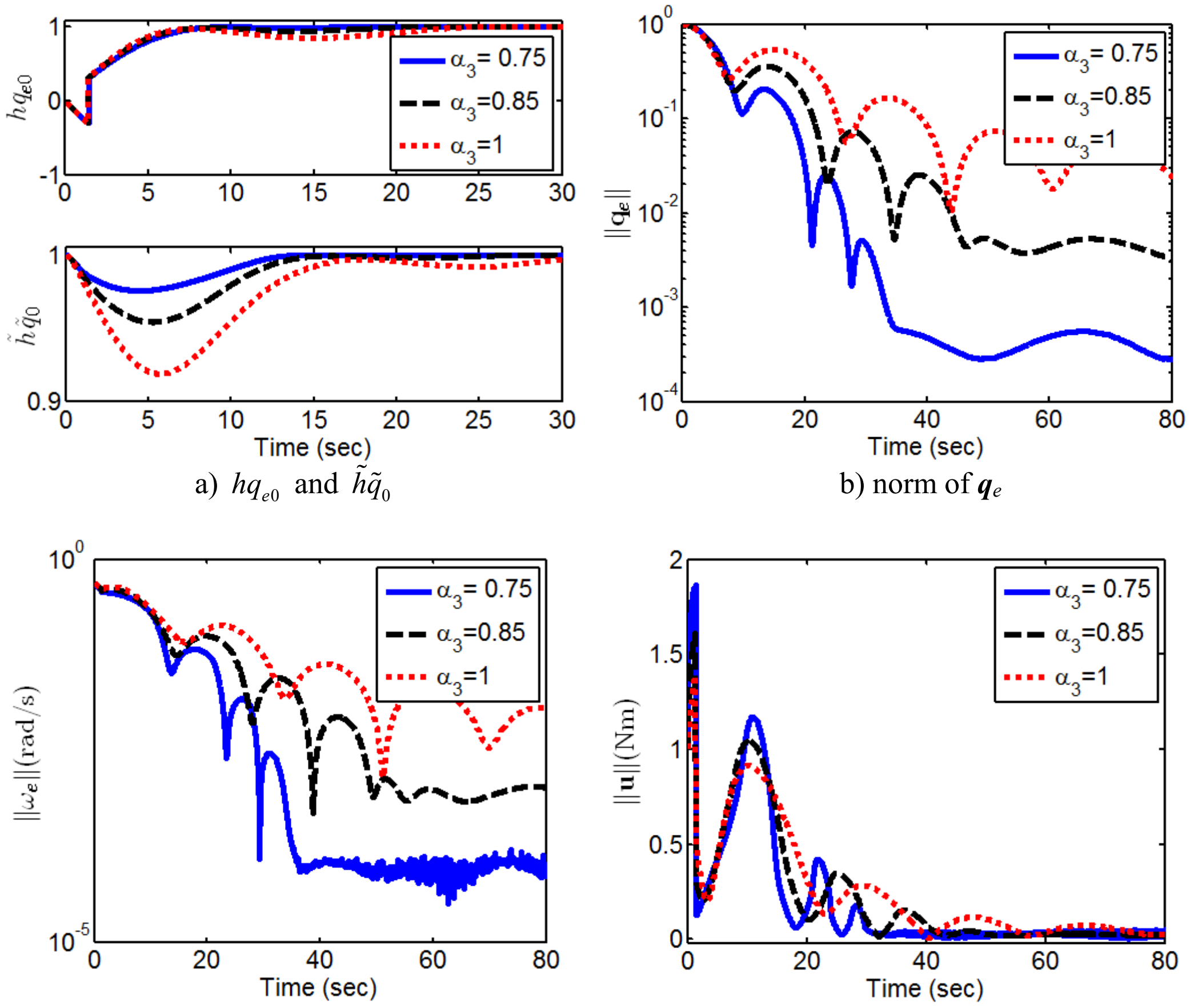

a) $hq_{e0}$ and $\tilde{h}\tilde{q}_0$ b) norm of $\boldsymbol{q}_e$

c) norm of $\boldsymbol{\omega}_e$ d) norm of $\boldsymbol{u}$

**Fig. 5 Simulation results for the GFTAC with attitude-only measurements and uncertainties.**

Figure 5 shows the simulation results for $\alpha_3 = 0.75, 0.85, 1$ ( $\alpha_1 = 2\alpha_3 - 1$ ). Similarly to the previous two examples, one jump occurs for $h$ while $\tilde{h}$ remains the same sign (Fig. 5a). As shown by the transient response of $\tilde{h}\tilde{q}_0$, the damping is successfully injected into the closed-loop system by the quaternion filter (36). Despite the absence of angular velocity measurements, $(\boldsymbol{Q}_e, \boldsymbol{\omega}_e)$ is stabilized to a small neighborhood of $(-\mathbf{1}, 0)$ and the tracking error reduces as the value of $\alpha_3$ decreases (Figs 5b and 5c). Clearly, the proposed GFATC achieves better robustness against system uncertainties than the asymptotic method in [7], which corresponds to $\alpha_1 = \alpha_3 = 1$.

## 5 Conclusions

The attitude tracking of a rigid body is studied in quaternion coordinates with full-state or attitude only measurements, or attitude plus biased angular velocity measurements. By dealing with these scenarios, we have shown how to combine the hybrid control method and homogeneous theory to achieve global finite-time attitude tracking while avoiding the noise-induced chattering of discontinuous switching laws. The resulting control laws all have a proportional-derivative form and ensure bounded control torques a priori. In particular, several existing attitude control schemes can be derived as special cases of the methods in this paper.

## Appendix A: Proof of Lemma 2.1

According to Definition 2.2, we only need to prove the uniform local finite-time stability of system (2) in order to obtain the result of Lemma 2.1. The proof is similar to that of Theorem 3 in [22] and is divided into two steps. The first step obtains a strict Lyapunov function $V(\boldsymbol{x})$ for the lower order system $\dot{\boldsymbol{x}} = \boldsymbol{f}(\boldsymbol{x})$ while the second shows that $V(\boldsymbol{x})$ is a locally strict Lyapunov function for system (2). The details are as follows.

Since $\dot{\boldsymbol{x}} = \boldsymbol{f}(\boldsymbol{x})$ is asymptotically stable and $\boldsymbol{f}(\boldsymbol{x})$ is homogeneous of degree $k < 0$ with respect to a dilation $\Delta_\varepsilon^r$, it follows from Theorem 7.1 in [27] that $\dot{\boldsymbol{x}} = \boldsymbol{f}(\boldsymbol{x})$ is finite-time stable. According to Theorem 7.2 in [27], there exist a positive-definite $C^1$ function $V : \mathbb{R}^n \to \mathbb{R}_{\geq 0}$ with homogeneous degree of $l > -k$, and a positive constant $c = -\max_{\{z:V(z)=1\}} L_f V(z) > 0$ such that

$$L_f V(\boldsymbol{x}) \leq -cV^\beta(\boldsymbol{x}),\ \boldsymbol{x} \in \mathbb{R}^n,\ \beta = (l+k)/l. \tag{A.1}$$

Given a constant $c_1$ satisfying $0 < c_1 < c$, it follows from (A.1) that

$$L_f V(\boldsymbol{x}) + c_1 V^\beta(\boldsymbol{x}) \le -(c - c_1) V^\beta(\boldsymbol{x}),$$

which implies $L_f V(\boldsymbol{x}) + c_1 V^\beta(\boldsymbol{x}) < 0$ for any $\boldsymbol{x} \ne 0$. Now, we can define a negative constant as $c_0 = \max\{L_f V(\boldsymbol{x}) + c_1 V^\beta(\boldsymbol{x}) : \boldsymbol{x} \in \mathbb{S}^{n-1}\} < 0$, where $\mathbb{S}^{n-1} = \{\boldsymbol{x} \in \mathbb{R}^n : \|\boldsymbol{x}\| = 1\}$ is a unit sphere on $\mathbb{R}^n$. In addition, the homogeneity of $V(\boldsymbol{x})$ and $\boldsymbol{f}(\boldsymbol{x})$ implies that $\partial V / \partial x_i$ and $L_f V(\boldsymbol{x})$ are also homogeneous since $\partial V(\Delta_\varepsilon^r \boldsymbol{x}) / \partial x_i = \varepsilon^{l - r_i} \partial V(\boldsymbol{x}) / \partial x_i$ and $L_f V(\Delta_\varepsilon^r \boldsymbol{x}) = \varepsilon^{l+k} L_f V(\boldsymbol{x})$ [22]. Applying this homogeneous property yields

$$L_{\hat{f}} V(\Delta_\varepsilon^r \boldsymbol{x}, t) = \sum_{i=1}^{n} \varepsilon^{l - r_i} \frac{\partial V(\boldsymbol{x})}{\partial x_i} \hat{f}_i(\Delta_\varepsilon^r \boldsymbol{x}, t) = \varepsilon^{l+k} \sum_{i=1}^{n} \frac{\partial V(\boldsymbol{x})}{\partial x_i} \frac{\hat{f}_i(\Delta_\varepsilon^r \boldsymbol{x}, t)}{\varepsilon^{r_i + k}}.$$

Note that $L_{\hat{f}} V$ depends explicitly on time since $\hat{\boldsymbol{f}}$ is time-dependent here. Invoking the condition given in (3) yields that $\lim_{\varepsilon \to 0} L_{\hat{f}} V(\Delta_\varepsilon^r \boldsymbol{x}, t) / \varepsilon^{l+k} = 0$ and thus $L_{\hat{f}} V(\Delta_\varepsilon^r \boldsymbol{x}, t) = o(\varepsilon^{l+k})$ uniformly holds on $\boldsymbol{x} \in \mathbb{S}^{n-1}$ for all $t \ge 0$. Hence, there exists $\varepsilon_0 \in (0,1)$ such that

$$\left| L_{\hat{f}} V(\Delta_\varepsilon^r \boldsymbol{x}, t) \right| \le -c_0 \varepsilon^{l+k} / 2, \tag{A.2}$$

holds for all $0 < \varepsilon < \varepsilon_0$, $\boldsymbol{x} \in \mathbb{S}^{n-1}$, and $t \ge 0$. Next, we compute the time derivative of $V$ along system (2), i.e., $\dot{V}$, which depends explicitly on time. It follows from (A.2) that

$$\begin{aligned} \dot{V}(\Delta_\varepsilon^r \boldsymbol{x}, t) + c_1 V^\beta(\Delta_\varepsilon^r \boldsymbol{x}) &= L_f V(\Delta_\varepsilon^r \boldsymbol{x}) + L_{\hat{f}} V(\Delta_\varepsilon^r \boldsymbol{x}, t) + c_1 V^\beta(\Delta_\varepsilon^r \boldsymbol{x}) \\ &= \varepsilon^{l+k} L_f V(\boldsymbol{x}) + c_1 \varepsilon^{l+k} V^\beta(\boldsymbol{x}) + L_{\hat{f}} V(\Delta_\varepsilon^r \boldsymbol{x}, t) \\ &\le \varepsilon^{l+k} [L_f V(\boldsymbol{x}) + c_1 V^\beta(\boldsymbol{x}) - 0.5 c_0] \le 0.5 \varepsilon^{l+k} c_0 < 0 \end{aligned}, \tag{A.3}$$

holds for all $0 < \varepsilon < \varepsilon_0$, $\boldsymbol{x} \in \mathbb{S}^{n-1}$, and $t \ge 0$. Define a set $U = \{\boldsymbol{z} : \boldsymbol{z} = \Delta_\varepsilon^r \boldsymbol{x}, \boldsymbol{x} \in \mathbb{S}^{n-1}, \varepsilon \in [0, \varepsilon_0)\}$, which can be obtained from $\Delta_\varepsilon^r \boldsymbol{x}$ ($\boldsymbol{x} \in \mathbb{S}^{n-1}$) by continuously varying $\varepsilon$ in $[0, \varepsilon_0)$. Clearly, $U \subseteq \mathbb{R}^n$ is a neighborhood of $\boldsymbol{x} = 0$ [22] and is inside the unit sphere $\mathbb{S}^{n-1}$ since $\varepsilon < \varepsilon_0 < 1$. Equation (A.3) then implies that $\dot{V}(\boldsymbol{x}, t) \le -c_1 V^\beta(\boldsymbol{x}) \le 0$ holds for all $(\boldsymbol{x}, t) \in U \times \mathbb{R}_{\ge 0}$. Hence, $\boldsymbol{x} = 0$ is uniformly asymptotically stable in $U$. Denote $\dot{V}(t) \triangleq \dot{V}(\boldsymbol{x}(t), t)$, $V(t) \triangleq V(\boldsymbol{x}(t))$, and $\boldsymbol{x}_0 \triangleq \boldsymbol{x}(t_0)$. Applying the comparison principle [23] to $\dot{V}(t) \le -c_1 V^\beta(t)$ and recognizing $0 < \beta < 1$ leads to $V^{1-\beta}(t) \le V^{1-\beta}(t_0) - c_1(1-\beta)(t - t_0)$. Hence, $V(t) = 0$ and $\boldsymbol{x}(t, \boldsymbol{x}_0, t_0) = 0$ for all $t - t_0 \ge V^{1-\beta}(t_0) / c_1(1-\beta)$. Noting $V(\boldsymbol{x}_0) = V(t_0)$, the settling time for any $(\boldsymbol{x}_0, t_0) \in U \times \mathbb{R}_{\ge 0}$ is then uniformly bounded by $T(\boldsymbol{x}_0) \le V^{1-\beta}(\boldsymbol{x}_0) / c_1(1-\beta)$. Therefore, system (2) is uniformly locally finite-time stable.

## Appendix B

**Proof of Property 1**: If $\|\boldsymbol{q}\| \ne 0$, (5) implies that $\boldsymbol{\kappa}_0(\boldsymbol{Q}, \alpha)$ is continuous and $\|\boldsymbol{\kappa}_0(\boldsymbol{Q}, \alpha)\| = \|\boldsymbol{q}\|^{1-\alpha}$. It then follows that $\lim_{\|\boldsymbol{q}\| \to 0} \|\boldsymbol{\kappa}_0(\boldsymbol{Q}, \alpha)\| = 0$ and thus $\boldsymbol{\kappa}_0(\boldsymbol{Q}, \alpha)$ is also continuous at $\|\boldsymbol{q}\| = 0$. In addition, $q_0^2 + \|\boldsymbol{q}\|^2 = 1$ implies that

$$2(1 - q_0) = (q_0 - 1)^2 + \|\boldsymbol{q}\|^2, \tag{B.1}$$

and thus $2(1 - q_0) \ge \|\boldsymbol{q}\|^2$. One can then deduce that $\|\boldsymbol{\kappa}_1(\boldsymbol{Q}, \alpha)\| \le \|\boldsymbol{q}\|^{1-\alpha} \le 1$ and $\lim_{q_0 \to 1} \|\boldsymbol{\kappa}_1(\boldsymbol{Q}, \alpha)\| = 0$. Hence, $\boldsymbol{\kappa}_1(\boldsymbol{Q}, \alpha)$ is continuous at the entire $\mathbb{S}^3$, too.

**Proof of Property 2**: It follows from (5) and (B.1) that $\boldsymbol{\kappa}_1(\boldsymbol{Q}, \alpha) = \boldsymbol{\kappa}_0(\boldsymbol{Q}, \alpha)(1 + x)^{-\alpha/2}$, where $x \triangleq (1 - q_0)^2 / \|\boldsymbol{q}\|^2$. If $q_0 > 0$, one can deduce that $1 - q_0 \le 1 - q_0^2 = \|\boldsymbol{q}\|^2$, which together with $0 \le 1 - q_0 \ll 1$ leads to $0 \le x \le 1 - q_0 \ll 1$. By a Taylor series expansion, it follows that

$$\boldsymbol{\kappa}_1(\boldsymbol{Q}, \alpha) - \boldsymbol{\kappa}_0(\boldsymbol{Q}, \alpha) \approx \boldsymbol{\kappa}_0(\boldsymbol{Q}, \alpha)(1 - \frac{\alpha}{2} x - 1) = -\frac{\alpha}{2} x \boldsymbol{\kappa}_0(\boldsymbol{Q}, \alpha), \tag{B.2}$$

and $1 - q_0 = 1 - \sqrt{1 - \|\boldsymbol{q}\|^2} \approx \|\boldsymbol{q}\|^2 / 2$. Property 2 is then verified by substituting $x \approx \|\boldsymbol{q}\|^2 / 4$ into (B.2).